\documentclass[10pt]{amsart}
\usepackage{latexsym, amsmath, amssymb}
\usepackage{epsfig,graphics}

\newcommand{\varespilon}{\varepsilon}

\newcommand{\newsection}[1]
{\section{#1}\setcounter{theorem}{0} \setcounter{equation}{0}
\par\noindent}

\newtheorem{theorem}{Theorem}

\newtheorem{lemma}[theorem]{Lemma}
\newtheorem{corr}[theorem]{Corollary}

\newtheorem{proposition}[theorem]{Proposition}

\newcommand{\cd}{\, \cdot\, }
\newcommand{\supp}{\text{supp }}
\newcommand{\R}{{\mathbb R}}
\newcommand{\ext}{{\R^3\backslash\mathcal{K}}}

\newcommand{\bdy}{{\partial\mathcal{K}}}

\newcommand{\K}{{\mathcal{K}}}

\renewcommand{\ss}[1]{{\substack{#1}}}

\newcommand{\la}{\langle}
\newcommand{\ra}{\rangle}

\newcommand{\1}{{\mathbf 1}}

\begin{document}

\title[Elastic waves in exterior domains]
{Elastic waves in exterior domains \\ Part II: Global existence with a null structure}

\thanks{The authors were supported by the NSF}

\author{Jason Metcalfe}
\address{Department of Mathematics, University of California, Berkeley, CA  94720-3840}
\email{metcalfe@math.berkeley.edu}

\author{Becca Thomases}
\address{Courant Institute of Mathematical Sciences, New York University, New York, NY  10012}
\email{thomases@cims.nyu.edu}

\begin{abstract}
In this article, we prove that solutions to a problem in nonlinear elasticity corresponding
to small initial displacements exist globally
in the exterior of a nontrapping obstacle.  The medium is assumed to be homogeneous, isotropic,
and hyperelastic, and the nonlinearity is assumed to satisfy a null condition.  The techniques contained
herein would allow for more complicated geometries provided that there is a sufficient decay of local energy
for the linearized problem.
\end{abstract}

\maketitle

\newsection{Introduction}

In this paper, we shall prove global existence of elastic waves in exterior domains subject to a null 
condition.  The elastic medium is assumed to be homogeneous, isotropic, and hyperelastic, and we 
only consider small initial displacements.  We solve the said equations exterior to a bounded, nontrapping
obstacle with smooth boundary, though our techniques would allow for any bounded, smooth obstacle for
which there is a sufficiently fast decay of local energy.  

For the boundaryless problem, global existence was previously shown in \cite{Agemi} and \cite{Si}.  
See, also, \cite{Si2, SiNotes}.  As in studies of the wave equation, the null condition is essential
in order to obtain global solutions.  Without some assumption on the structure of the nonlinearity,
blow up in finite time was shown in \cite{J2}.  In this case, almost global existence of solutions
was shown in \cite{J} and a simplified proof was later given in \cite{KS}.  In \cite{Met1}, the corresponding
almost global existence result was shown outside of a star-shaped obstacle.\footnote{By finite propagation
speed and \cite{J2}, blow up in finite time is still possible.}

Several related studies have been conducted concerning multiple speed systems of wave equations 
in three-dimensional exterior domains.  Almost global existence was established in
\cite{KSS2, KSS3} and \cite{MS3} exterior to star-shaped domains.  Using the techniques of
\cite{MS1}, these results can be extended to any domain for which there is a sufficiently fast
decay of local energy.  For nonlinearities satisfying the null condition, global existence has been
established in \cite{MS1, MS4} and \cite{MNS1, MNS2}.  Our proof will use many of the techniques
established in this context.

We shall utilize the method of commuting vector fields, though the equation of interest is not
Lorentz invariant and thus we will be required to use a restricted set of vector fields.  Moreover,
we will rely upon the arguments of \cite{KSS3} and \cite{MS1} to establish the necessary higher-order
energy estimates.  Here, we use elliptic regularity to establish energy estimates involving 
the generators of translations.  Without the star-shapedness assumption, proving energy estimates
involving the scaling vector field is more delicate, and we will use a boundary term estimate which
is reminiscent of those in \cite{MS1}.  For energy estimates involving the full set of vector fields,
a class of weighted mixed-norm estimates, called KSS estimates, will be called upon to handle the relevant
boundary terms.  These estimates are analogous to those of \cite{KSS2} for the wave equation and
have previously be established in \cite{Met1} for the current setting.

The nontrapping hypothesis is assumed to guarantee the existence of an exponential decay of local energy.
See \cite{Y}.  Similar well-known estimates (see \cite{MRS}) were widely used in the studies of
nonlinear wave equations in exterior domains.  Our techniques would, in fact, permit any exterior domain\footnote{By 
this, we mean the exterior of any bounded obstacle with smooth boundary.} for which
there is a sufficiently fast decay of local energy.  In particular, the method of proof would allow
for the loss of regularity in the local energy decay which would be necessary if there were trapped
rays.  In the setting of the wave equation, see \cite{I1, I2} for such local energy decays and
\cite{MS1}, \cite{MNS1, MNS2} for proofs of long-time existence in such domains.  Moreover, by
using the technique of proof contained herein, the almost global existence result of \cite{Met1}
may be extended to any domain for which there is such a decay of local energy.

Let us more precisely describe the initial-boundary value problem.  We fix a bounded, nontrapping\footnote{Recall 
that more general domains are possible, but the authors are not aware of the necessary
decay of local energy results for such domains.} obstacle $\K\subset\R^3$ with smooth boundary.  We may,
without loss of generality, assume that $0\in\K\subset\{|x|<1\}$, and we will do so throughout.

The linearized equation of elasticity is
\begin{equation}\label{linear}
Lu=\partial_t^2 u -c_2^2 \Delta u - (c_1^2-c_2^2)\nabla (\nabla\cdot u)=0,
\end{equation}
where $u=(u^1,u^2,u^3)$ is the displacement vector.  The constants $c_1, c_2$ may be assumed
to satisfy 
$$c_1^2-\frac{4}{3}c_2^2>0,\quad c_2^2>0.$$
The two quantities in the above equations represent the bulk and shear moduli respectively.  

The nonlinear problem that results from hyperelastic mediums is
\begin{equation}
\label{nonlinear}
(Lu)^I=\partial_l (B^{IJK}_{lmn}\partial_m u^J \partial_n u^K)+\dots\footnotemark
\end{equation}
\footnotetext{Here and throughout, we use the summation convention.  Latin indices indicate that
the implicit summation runs from $1$ to $3$, while Greek indices are used when the summation
shall run from $0$ to $3$.}where the real constants $B^{IJK}_{lmn}$ satisfy the symmetry condition
\begin{equation}
\label{nonlinear.symmetry}
B^{IJK}_{lmn}=B^{JIK}_{mln}=B^{IKJ}_{lnm}.
\end{equation}
The terms not explicitly stated in \eqref{nonlinear} are all of cubic or higher order in $\nabla u,
\nabla^2 u$.  The null condition that we require assumes that
\begin{equation}
\label{null.condition}
B^{IJK}_{lmn}\xi_I\xi_J\xi_K\xi_l\xi_m\xi_n=0,\quad \text{whenever } \xi\in S^2.
\end{equation}
Rather than interpreting the equation and null condition further, we refer the reader to 
\cite{Agemi}, \cite{Si}, and the references therein.  The interested reader may also
wish to refer to the expository article \cite{SiNotes}.  

As has become common, for convenience, we will truncate the equations at the quadratic level.  
Such a truncation does not affect the long-time existence.  We then study the following
initial-boundary value problem:
\begin{equation}
\label{main.equation}
\begin{cases}
(Lu)^I = \partial_l (B^{IJK}_{lmn}\partial_m u^J \partial_n u^K),\quad (t,x)\in \R_+\times
\ext,\\
u|_\bdy=0,\\
u(0,\cd)=f,\quad \partial_t u(0,\cd)=g
\end{cases}
\end{equation}
for ``small'' Cauchy data $f,g$. 
For convenience in the sequel,
we shall set
$$Q^I(\nabla u,\nabla^2 u) = \partial_l (B^{IJK}_{lmn}\partial_m u^J \partial_n u^K).$$

In order to solve \eqref{main.equation}, we must assume that the initial data satisfy
certain well-known compatibility conditions.  Letting $J_ku=\{\partial_x^\alpha u\,:\,0\le|\alpha|\le k\}$
and noticing that $\partial_t^k u(0,\cd)=\psi_k(J_kf,J_{k-1}g)$, $0\le k\le m$ for formal $H^m$ solutions,
we require that the compatibility functions $\psi_k$ -- which depend on $Q$, $J_kf$, and $J_{k-1}g$ --
vanish on $\bdy$ for $0\le k\le m-1$.  Smooth data $(f,g)\in C^\infty$ are said to satisfy
the compatibility conditions to infinite order if this holds for all $m$.

Under these assumptions, we may prove that solutions to \eqref{main.equation} corresponding to small
initial displacements exist globally.  This is our main result.
\begin{theorem}
\label{main.theorem}  Let $\K$ be a bounded, nontrapping obstacle with smooth boundary as above.  Assume
that $B^{IJK}_{lmn}$ are real constants satisfying \eqref{nonlinear.symmetry} and \eqref{null.condition}.
Suppose further that the data $(f,g)\in C^\infty(\ext)$ satisfy the compatibility conditions
to infinite order.  Then, there are positive constants $\varepsilon_0$ and $N$ so that for
all $\varepsilon\le \varepsilon_0$, if 
\begin{equation}
\label{data.smallness}
\sum_{|\alpha|\le N} \|\la x\ra^{|\alpha|}\partial_x^\alpha f\|_2
+\sum_{|\alpha|\le N-1} \|\la x\ra^{|\alpha|+1} \partial_x^\alpha g\|_2\le \varepsilon,\footnotemark
\end{equation}
then \eqref{main.equation} has a unique global solution $u\in C^\infty([0,\infty)\times \ext)$.
\footnotetext{Here, and throughout, $\la x\ra=\la r\ra = \sqrt{1+|x|^2}$ denotes the Japanese bracket.}
\end{theorem}

This paper is organized as follows.  In the remainder of this section, we introduce some notations
that will be used throughout.  In particular, we introduce the vector fields that shall be utilized.
The second section is devoted to estimates related to the energy inequality.  These include
the exponential decay of local energy, the energy inequality and its higher order variants, and weighted
mixed-norm KSS estimates.  In the third section, we gather the main pointwise decay estimates
that are used to give global existence.  As a corollary to these, we obtain the main boundary term
estimate that is utilized instead of the star-shapedness assumption which was used, e.g., in \cite{Met1}.
The fourth section is devoted to certain Sobolev-type estimates and estimates related to the null condition.
Finally, in the last section, we prove Theorem \ref{main.theorem}.

\noindent{\em Acknowledgements:} The authors are particularly grateful to T. Sideris and C. Sogge
for introducing us to this problem and for numerous related discussions and collaborations.

\subsection{Notation}
When convenient, we shall set $x_0=t$, $\partial_0=\partial_t$.  The space-time gradient
will be denoted by $u'=\partial u=(\partial_t, \nabla_x u)$, but we shall reserve the notation
$\nabla$ for
$\nabla=\nabla_x
=(\partial_1,\partial_2,\partial_3)$.  The generators of the spatial rotations are denoted
$$\Omega=(\Omega_1,\Omega_2,\Omega_3)=x\times \nabla,$$
where $\times$ is the usual vector cross product. 

When studying elasticity, it is natural to use the generators of the simultaneous rotations
$$\tilde{\Omega}_l = \Omega_l I + U_l,$$
where
$$U_1=\begin{bmatrix}0&0&0\\0&0&1\\0&-1&0\end{bmatrix},\quad
U_2=\begin{bmatrix}0&0&-1\\0&0&0\\1&0&0\end{bmatrix},\quad
U_3=\begin{bmatrix}0&1&0\\-1&0&0\\0&0&0\end{bmatrix}.$$

The scaling operator will be denoted
$$S=t\partial_t + r\partial_r -1.$$
This differs slightly from the usual notion of the scaling vector field, but here it is more
convenient to use this modification as it properly preserves the null structure.  See \cite{Si}.

We will set
$$Z=\{\partial, \tilde{\Omega}\}, \quad \Gamma=\{\partial,\tilde{\Omega}, S\}.$$
We differentiate these two sets as it will be necessary to produce estimates which require
relatively few occurrences of the scaling vector field $S$.  This is due to the fact that
the coefficients of $S$ can be arbitrarily large in a neighborhood of $\bdy$, while those of $Z$
are bounded.  Similar care with $S$ had to be taken in \cite{KSS3}, \cite{MS1}, \cite{MNS1, MNS2},
and \cite{Met1}.

A key property of the vector fields is the commutation properties with $L$.  In particular, we have
that 
\begin{equation}\label{commutators}
[L,\partial]=[L,\tilde{\Omega}]=0,\quad [L,S]=2L.\end{equation}
Moreover, we note that the vector fields preserve the null structure of $Q$.  See \cite{Si}.\footnote{See
Proposition 3.1.}

We will use the projections onto the radial and transverse directions, defined as follows:
$$P_1 u = \frac{x}{r}\la \frac{x}{r}, u\ra,
\quad P_2 u = [I-P_1]u = -\frac{x}{r}\times\Bigl(\frac{x}{r}\times u\Bigr).$$

We shall use $A\lesssim B$ to denote that there is a positive, unspecified constant $C$ so that
$A\le CB$, and we will use $S_t=[0,t]\times\ext$ to denote a time-strip of height $t$ in the exterior 
domain.

\bigskip
\newsection{Energy and KSS estimates}

\subsection{Local energy decay}
An important tool in previous studies of nonlinear problems in exterior domains is the decay of local
energy.  It is here that we require the geometric condition on the obstacle.
In the current setting, we have the following result of \cite{Y}.  This is an analog of the classical 
result of \cite{MRS} for the wave equation.
\begin{theorem}\label{theorem.local.energy}
 Let $\K\subset\{|x|<1\}\subset \R^3$ be a nontrapping obstacle with smooth boundary.  Let
$u$ solve 
\begin{equation}
\label{local.energy.equation}
\begin{cases}
Lu=0,\quad (t,x)\in \R_+\times\ext,\\
u|_\bdy=0,\\
\supp u(0,\cd), \partial_t u(0,\cd)\subset \{|x|<10\}.
\end{cases}
\end{equation}
Then
\begin{equation}
\label{local.energy}
\Bigl(\int_{\{x\in\ext\,:\, |x|<10\}} |u'(t,x)|^2\:dx\Bigr)^{1/2}\lesssim e^{-ct} \|u'(0,\cd)\|_2,
\end{equation}
for some $c>0$.\end{theorem}

In what follows, we shall require a higher order version of \eqref{local.energy}.  To establish this,
we utilize the following version of elliptic regularity, which appeared in \cite{Met1}.
\begin{lemma}
\label{lemma.ell.reg}
Let $\K\subset\{|x|<1\}\subset \R^3$ be an obstacle with smooth boundary.  Suppose that $u\in C^\infty
(\R_+\times\ext)$, $u|_\bdy=0$, and $u$ vanishes for large $|x|$ for each $t$.  Then,
\begin{equation}
\label{ell.reg}
\sum_{|\alpha|\le M}\|S^\nu \partial^\alpha u'(t,\cd)\|_2\lesssim
\sum_\ss{j+\mu\le M+\nu\\\mu\le\nu} \|S^\mu \partial_t^j u'(t,\cd)\|_2
+\sum_\ss{|\beta|+\mu\le M+\nu-1\\\mu\le\nu} \|S^\mu \partial^\beta Lu(t,\cd)\|_2
\end{equation}
and
\begin{multline}
\label{ell.reg.local}
\sum_{|\alpha|\le M} \|S^\nu \partial^\alpha u'(t,\cd)\|_{L^2(|x|<4)} \lesssim
\sum_\ss{j+\mu\le M+\nu\\\mu\le\nu} \|S^\mu \partial_t^j u'(t,\cd)\|_{L^2(|x|<6)}
\\+ \sum_\ss{|\beta|+\mu\le M+\nu-1\\\mu\le\nu} \|S^\mu \partial^\beta Lu(t,\cd)\|_{L^2(|x|<6)}
\end{multline}
for any $M$ and $\nu$.
\end{lemma}

Using \eqref{ell.reg.local} and the fact that $\partial_t$ preserves the Dirichlet boundary conditions,
we can immediately establish the following result which is also from \cite{Met1}.
\begin{lemma}
\label{lemma.local.energy.high} 
Let $\K\subset \{|x|<1\}\subset \R^3$ be a smooth, nontrapping obstacle.
Suppose that $u|_\bdy=0$ and  $Lu(t,x)=0$ for $|x|>4$ and $t>0$.  
Suppose also that $u(t,x)=0$ for $t\le 0$.  Then if $M$ and $\nu$ are fixed and
if $c>0$ is as in \eqref{local.energy},
\begin{multline}
\label{local.energy.high}
\sum_\ss{|\alpha|+\mu\le M+\nu\\\mu\le \nu} \|S^\mu \partial^\alpha u'(t,\cd)\|_{L^2(|x|<4)}
\lesssim 
\sum_\ss{|\alpha|+\mu\le M+\nu-1\\\mu\le\nu} \|S^\mu \partial^\alpha Lu(t,\cd)\|_2
\\+\int_0^t e^{-(c/2)(t-s)} \sum_\ss{|\alpha|+\mu\le M+\nu\\\mu\le\nu} \|S^\mu \partial^\alpha
Lu(s,\cd)\|_2\:ds.
\end{multline}
\end{lemma}

\subsection{Energy estimates}\label{energy.estimates.section}
In this section, we gather the energy estimates which we shall require.  The basic energy
estimate is rather standard.  In order to establish estimates for $\partial^\alpha u$, we shall use
elliptic regularity as above.  For energy estimates involving the scaling vector field, we use
cutoff techniques from \cite{MS1} and control the resulting commutator with estimates given in Section
\ref{boundary}.  Finally, we handle the boundary terms that arise when studying $S^\mu Z^\alpha u$
using the KSS estimates given in the following section.

We begin with the standard energy inequality for the variable coefficient operator
\begin{equation}
\label{Lgamma}
(L_\gamma u)^I = (Lu)^I + \gamma^{IJ,jk}\partial_j\partial_k u^J.
\end{equation}
We look at smooth solutions of
\begin{equation}
\label{perturbed.equation}
\begin{cases}
L_\gamma u = F,\quad (t,x)\in \R_+\times\ext,\\
u|_\bdy=0,\\
u(0,\cd)=f,\quad \partial_tu(0,\cd)=g
\end{cases}
\end{equation}
assuming that 
\begin{equation}
\label{gamma.symmetry}
\gamma^{IJ,jk}=\gamma^{JI,kj}
\end{equation}
and 
\begin{equation}
\label{gamma.smallness}
\sum_{I,J=1}^3 \sum_{j,k=1}^3 \|\gamma^{IJ,jk}(t,\cd)\|_\infty\le \frac{\delta}{1+t}
\end{equation}
for $\delta>0$ sufficiently small, depending on $c_1$ and $c_2$.

We define the energy-momentum vector associated to $L_\gamma$,
\begin{equation}
\label{e0}
e_0[u]=\frac{1}{2}|\partial_t u|^2+\frac{c_2^2}{2}|\nabla u|^2 + \frac{c_1^2-c_2^2}{2}(\nabla\cdot u)^2
-\frac{1}{2}\gamma^{IJ,ij}\partial_i u^I \partial_j u^J,
\end{equation}
\begin{equation}\label{ek}
e_k[u]=-c_2^2\partial_t u^I \partial_k u^I - (c_1^2-c_2^2)\partial_t u^k (\nabla\cdot u)
+\gamma^{IJ,kj}\partial_j u^J \partial_t u^I,\quad k=1,2,3.
\end{equation}
With \eqref{gamma.symmetry}, it is easy to check that
\begin{equation}
\label{e.div}
\partial_0 e_0[u]+\partial_k e_k[u]=\partial_t u^I (L_\gamma u)^I
-\frac{1}{2}(\partial_t \gamma^{IJ,ij})\partial_i u^I \partial_j u^J
+ (\partial_k \gamma^{IJ,kj})\partial_j u^J \partial_t u^I.
\end{equation}

Setting
$$E_M(t)=\int \sum_{j=0}^M e_0(\partial_t^j u)(t,x)\:dx,$$
the following results immediately from \eqref{e.div} and that $\partial_t$ preserves the boundary conditions.
\begin{lemma}
\label{lemma.energy.dt}
Assume that $\gamma$ satisfies \eqref{gamma.symmetry} and \eqref{gamma.smallness}, and let $u$ be a smooth solution to 
\eqref{perturbed.equation} which vanishes for large $|x|$ for each $t$.  Then
\begin{equation}
\label{energy.dt}
\partial_t E^{1/2}_M(t)\lesssim \sum_{j=0}^M \|L_\gamma \partial_t^j u(t,\cd)\|_2
+ \|\gamma'(t,\cd)\|_\infty E_M^{1/2}(t)
\end{equation}
for any fixed $M=0,1,2,\dots$.\footnote{Here, we have set 
$\|\gamma'(t,\cd)\|_\infty = \sum_{I,J=1}^3\sum_{j,k=0}^3 \sum_{\beta=0}^3 \|\partial_\beta
\gamma^{IJ,jk}(t,\cd)\|_\infty.$}
\end{lemma}

We next examine energy estimates involving the scaling vector field.  To do so, we set
$$\tilde{S}=t\partial_t + \eta(x) r\partial_r -1$$
where $\eta\in C^\infty(\R^3)$ with $\eta(x)\equiv 0$ for $x\in \K$ and $\eta(x)\equiv 1$ for $|x|>1$.
For
$$X_{\nu,j}=\int e_0(\tilde{S}^\nu \partial_t^j u)(t,x)\:dx,$$
we have the following.
\begin{lemma}
\label{lemma.energy.tilde.S}
If $u$ is a smooth solution to \eqref{perturbed.equation} and vanishes for large $|x|$ for each $t$,
then
\begin{multline}
\label{energy.tilde.S}
\partial_t X_{\nu,j}\lesssim X^{1/2}_{\nu,j} \|\tilde{S}^\nu \partial_t^j L_\gamma u(t,\cd)\|_2
+ \|\gamma'(t,\cd)\|_\infty X_{\nu,j}
\\+X_{\nu,j}^{1/2} \|[\tilde{S}^\nu \partial_t^j,\gamma^{kj}\partial_k\partial_l]u(t,\cd)\|_2
+ X_{\nu,j}^{1/2}\sum_{\mu\le\nu-1} \|S^\mu \partial_t^j Lu(t,\cd)\|_2
\\+ X^{1/2}_{\nu,j}\sum_\ss{|\alpha|+\mu\le j+\nu\\\mu\le\nu-1} \|S^\mu \partial^\alpha u'(t,\cd)
\|_{L^2(\{|x|<1\})}
\end{multline}
for any fixed $\nu,j$.\footnote{For a differential operator $P=P(t,x,D_t,D_x)$, we set
$[P,\gamma^{kl}\partial_k\partial_l]u=\sum_{I,J=1}^3 \sum_{k,l=1}^3 |[P,\gamma^{IJ,kl}\partial_k\partial_l]
u^J|.$}
\end{lemma}

\noindent{\em Proof of Lemma \ref{lemma.energy.tilde.S}:}
As $\tilde{S}$ and $\partial_t$ preserve the boundary condition, it follows from \eqref{energy.dt} that
\begin{equation}
\label{energy.tilde.S.1}
\partial_t X_{\nu,j}\lesssim X_{\nu,j}^{1/2} \|L_\gamma \tilde{S}^\nu \partial_t^j u(t,\cd)\|_2
+ \|\gamma'(t,\cd)\|_\infty X_{\nu,j}.
\end{equation}
We next examine the commutator
\begin{align*}
|L_\gamma \tilde{S}^\nu \partial_t^j u| &\le |\tilde{S}^\nu \partial_t^j L_\gamma u|
+ |[\tilde{S}^\nu \partial_t^j, \gamma^{kl}\partial_k\partial_l]u|
+ |[\tilde{S}^\nu, L]\partial_t^j u|\\
&\le |\tilde{S}^\nu \partial_t^j L_\gamma u|+ |[\tilde{S}^\nu \partial_t^j, \gamma^{kl}\partial_k
\partial_l] u| + |[S^\nu,L]\partial_t^j u| + |[\tilde{S}^\nu - S^\nu, L]\partial_t^j u|\\
&\lesssim |\tilde{S}^\nu \partial_t^j L_\gamma u| + |[\tilde{S}^\nu \partial_t^j, \gamma^{kl}
\partial_k\partial_l]u|+\sum_{\mu\le\nu-1} |S^\mu \partial_t^j Lu|\\
&\qquad\qquad\qquad\qquad\qquad\qquad\qquad\qquad + \1_{\{|x|<1\}}(x) 
\sum_\ss{|\alpha|+\mu\le j+\nu\\\mu\le\nu-1} |S^\mu \partial^\alpha u'|.
\end{align*}
Using this in \eqref{energy.tilde.S.1} completes the proof.  \qed

Using the previous lemma and elliptic regularity, we shall prove
\begin{proposition}
\label{prop.energy.L.d}
Assume \eqref{gamma.symmetry} and \eqref{gamma.smallness} for $\delta$ sufficiently small, and let 
$u$ be a smooth solution to \eqref{perturbed.equation}.  Suppose further that
\begin{equation}
\label{gamma.smallness.2}
\|\gamma'(t,\cd)\|_\infty\le \delta/(1+t)
\end{equation}
and
\begin{multline}
\label{dt.bound}
\sum_\ss{j+\mu\le N+\nu\\\mu\le\nu} \Bigl(\|\tilde{S}^\mu \partial_t^j L_\gamma u(t,\cd)\|_2
+ \|[\tilde{S}^\mu \partial_t^j, \gamma^{kl}\partial_k\partial_l]u(t,\cd)\|_2\Bigr)
\\\le \frac{\delta}{1+t} \sum_\ss{j+\mu\le N+\nu\\\mu\le\nu} \|\tilde{S}^\mu \partial_t^j u'(t,\cd)\|_2
+ H_{\nu,N}(t)
\end{multline}
for some fixed $N$ and $\nu$ and some function $H_{\nu,N}(t)$.\footnote{In the sequel, $H_{\nu,N}(s)$ will
involve $\|\la x\ra^{-1/2} S^\mu Z^\alpha u'(s,\cd)\|^2_2$ for $|\alpha|+\mu\ll N+\nu$.  As such, this term
will be handled using the KSS estimates in the next section.}  Then,
\begin{multline}
\label{energy.L.d}
\sum_\ss{|\alpha|+\mu\le N+\nu\\\mu\le\nu} \|S^\mu \partial^\alpha u'(t,\cd)\|_2
\lesssim \sum_\ss{|\alpha|+\mu\le N+\nu-1\\\mu\le\nu} \|S^\mu \partial^\alpha Lu(t,\cd)\|_2
\\+ (1+t)^{A\delta} \sum_\ss{j+\mu\le N+\nu\\\mu\le\nu} X^{1/2}_{\mu,j}(0)
\\+(1+t)^{A\delta}\Bigl(\int_0^t \sum_\ss{|\alpha|+\mu\le N+\nu-1\\\mu\le \nu-1} \|S^\mu \partial^\alpha
Lu(s,\cd)\|_2\:ds + \int_0^t H_{\nu,N}(s)\:ds\Bigr)
\\+(1+t)^{A\delta} \int_0^t \sum_\ss{|\alpha|+\mu\le N+\nu\\\mu\le\nu-1} \|S^\mu \partial^\alpha u'
(s,\cd)\|_{L^2(|x|<1)}\:ds,
\end{multline}
for some constant $A>0$.
\end{proposition}

\noindent{\em Proof of Proposition \ref{prop.energy.L.d}:}
For $\delta$ in \eqref{gamma.smallness} sufficiently small,
\begin{equation}\label{energy.sim}
\sum_\ss{j+\mu\le N+\nu\\\mu\le \nu} \|\tilde{S}^\mu \partial_t^j u'(t,\cd)\|_2 \lesssim
\sum_\ss{j+\mu\le N+\nu\\\mu\le\nu} X^{1/2}_{\mu,j}(t).\end{equation}
By \eqref{energy.tilde.S}, \eqref{gamma.smallness.2}, and \eqref{dt.bound}, it thus follows that
\begin{multline}
\label{energy.L.d.1}
\partial_t \sum_\ss{j+\mu\le N+\nu\\\mu\le\nu} X^{1/2}_{\mu,j}(t)\le \frac{A\delta}{1+t}
\sum_\ss{j+\mu\le N+\nu\\\mu\le\nu} X^{1/2}_{\mu,j}(t)+AH_{\nu,N}(t) 
\\+A\sum_\ss{j+\mu\le N+\nu-1\\\mu\le\nu-1} \|S^\mu \partial_t^j Lu(t,\cd)\|_2
+A\sum_\ss{|\alpha|+\mu\le N+\nu\\\mu\le\nu-1} \|S^\mu \partial^\alpha u'(t,\cd)\|_{L^2(|x|<1)}.
\end{multline}
By Gronwall's inequality, $\sum_\ss{j+\mu\le N+\nu\\\mu\le\nu} X^{1/2}_{\mu,j}$ satisfies the desired
bound. 
Thus, by combining this estimate for $X^{1/2}_{\mu,j}$, \eqref{energy.sim}, \eqref{ell.reg}, 
and \eqref{energy.L.d.1}, the estimate \eqref{energy.L.d} follows.\qed

Finally, we establish the necessary energy estimates for $S^\mu Z^\alpha u$.  To this end, we set
$$Y_{N,\nu}(t)=\sum_\ss{|\alpha|+\mu\le N+\nu\\\mu\le\nu} \int e_0(S^\mu Z^\alpha u)(t,x)\:dx.$$
We then argue as in the proof
of Lemma \ref{lemma.energy.dt} with $E_M(t)$ replaced by $Y_{N,\nu}(t)$.  
The boundary terms that arise mesh well with the KSS estimates of the next section.
\begin{proposition}
\label{prop.energy.S.Z}  Assume \eqref{gamma.symmetry} and \eqref{gamma.smallness}.
Suppose that $u$ solves \eqref{perturbed.equation} and vanishes for large $|x|$ for each $t$.  Then,
\begin{multline}
\label{energy.S.Z}
\partial_t Y_{N,\nu}\lesssim Y_{N,\nu}^{1/2} \sum_\ss{|\alpha|+\mu\le N+\nu\\\mu\le\nu}
\|L_\gamma S^\mu Z^\alpha u(t,\cd)\|_2
+ \|\gamma'(t,\cd)\|_\infty Y_{N,\nu}\\+ \sum_\ss{|\alpha|+\mu\le N+\nu+1\\\mu\le\nu}
\|S^\mu \partial^\alpha u'(t,\cd)\|^2_{L^2(|x|<1)}.
\end{multline}
\end{proposition}

\noindent{\em Proof of Proposition \ref{prop.energy.S.Z}:}
By arguing as in the proof of Lemma \ref{lemma.energy.dt}, it follows that
\begin{multline}
\label{Y.div}
\partial_t Y_{N,\nu} \lesssim Y_{N,\nu}^{1/2} \sum_\ss{|\alpha|+\mu\le N+\nu\\\mu\le\nu}
\|L_\gamma S^\mu Z^\alpha u(t,\cd)\|_2 + \|\gamma'(t,\cd)\|_\infty Y_{N,\nu}
\\+ \sum_\ss{|\alpha|+\mu\le N_0+\nu_0\\\mu\le\nu_0} \int_\bdy |e_k(S^\mu Z^\alpha u)(t,y) n_k|\:d\sigma(y)
\end{multline}
where $n=(n_1,n_2,n_3)$ is the outward normal to $\K$ at a point $x\in \bdy$ and $e_k[\cd]$ are as
in \eqref{ek}.

Recalling that $\K\subset \{|x|<1\}$, 
$$\sum_\ss{|\alpha|+\mu\le N+\nu\\\mu\le\nu} |S^\mu Z^\alpha u(t,x)|\lesssim
\sum_\ss{|\alpha|+\mu\le N+\mu\\\mu\le\nu} |S^\mu \partial^\alpha u(t,x)|,\quad x\in\bdy,$$
and thus, by a trace theorem, we have that the last term in \eqref{Y.div} is
$$\lesssim \int_{\{x\in\ext\,:\,|x|<1\}} \sum_\ss{|\alpha|+\mu\le N+\nu+1\\\mu\le \nu}
|S^\mu \partial^\alpha u'(t,x)|^2\:dx,$$
which completes the proof.\qed

\subsection{KSS estimates}
In this section, we present a class of weighted, mixed-norm estimates, called KSS estimates, which
are particularly useful for estimating the boundary term in \eqref{energy.S.Z} and for dealing
with certain technicalities regarding the distribution of the occurrences of the scaling vector field
in the proof of the main theorem.  Such estimates were first used to study nonlinear problems
in \cite{KSS2} and have played a fundamental role in previous studies of problems in exterior domains.

The estimates that we give are from \cite{Met1}, and we refer the interested reader to that article for
detailed proofs.  They are based on the boundaryless estimates
for the wave equation from \cite{Sterb} and \cite{MS3}.  The corresponding estimates for elasticity
in the boundaryless case follows from a Helmholtz-Hodge decomposition.  In order to prove estimates
in the exterior domain one relies on the decay of local energy when $x$ is near the boundary
of the obstacle and uses the boundaryless estimates when $|x|$ is large.

\begin{proposition}
\label{prop.kss}
Suppose that $\K\subset\{|x|<1\}\subset\R^3$ is a nontrapping obstacle with smooth boundary.
Suppose further that $u\in C^\infty$ satisfies $u|_\bdy=0$, $u(t,x)=0$ for $t\le 0$, and vanishes
for large $|x|$ for every $t$.  Then,
\begin{multline}
\label{KSS.S.d}
(\log(2+T))^{-1/2} \sum_\ss{|\alpha|+\mu\le M+\nu\\\mu\le\nu} \|\la x\ra^{-1/2} S^\mu \partial^\alpha
u'\|_{L^2_tL^2_x(S_T)} \\\lesssim
\int_0^T \sum_\ss{|\alpha|+\mu\le M+\nu\\\mu\le\nu} \|S^\mu \partial^\alpha Lu(s,\cd)\|_2\:ds
+\sum_\ss{|\alpha|+\mu\le M+\nu-1\\\mu\le\nu} \|S^\mu \partial^\alpha Lu\|_{L^2_tL^2_x(S_T)}
\end{multline}
and
\begin{multline}
\label{KSS.S.Z}
(\log(2+T))^{-1/2}\sum_\ss{|\alpha|+\mu\le M+\nu\\\mu\le \nu} \|\la x\ra^{-1/2} S^\mu Z^\alpha
u'\|_{L^2_tL^2_x(S_T)}\\\lesssim
\int_0^T \sum_\ss{|\alpha|+\mu\le M+\nu\\\mu\le\nu} \|S^\mu Z^\alpha Lu(s,\cd)\|_2\:ds
+\sum_\ss{|\alpha|+\mu\le M+\nu-1\\\mu\le\nu} \|S^\mu Z^\alpha Lu\|_{L^2_tL^2_x(S_T)}
\end{multline}
for any fixed $M,\nu$ and $T>0$.
\end{proposition}

\bigskip
\newsection{Pointwise estimates and boundary term estimates}

In this section, we shall present our main pointwise decay estimates.  These are analogues of 
those used in \cite{KSS3} for the wave equation.  See also \cite{MS1}, \cite{MNS1, MNS2}.  In the process
of proving such estimates, we shall also prove the boundary term estimate which is required to handle the
last term of \eqref{energy.L.d}.  
This is reminiscent of ideas from \cite{MS1, MS2}.  

We begin by looking at the solution to the boundaryless problem
\begin{equation}
\label{boundaryless.equation}
\begin{cases}
Lv(t,x)=0,\\
v(0,\cd)=f,\quad \partial_t v(0,\cd)=g.
\end{cases}
\end{equation}
By arguing using spherical means\footnote{See \cite{J}.}, we have that
\begin{multline}
\label{fund.soln}
4\pi v_I(t,x)=\int_{S^2} \Bigl[(2\delta_{IJ} - 4 y_I y_J)f_J(x+c_2ty)
\\\qquad\qquad\qquad\qquad\qquad\qquad\qquad
+t (\delta_{IJ}-y_Iy_J)(g_J(x+c_2ty)+c_2 y_K \nabla_Kf_J(x+c_2ty))\Bigr]\:d\sigma(y)
\\\!\!\!\!\!\!\!\!\!\!\!\!\!\!\!\!+\int_{S^2} \Bigl[(-\delta_{IJ}+4y_Iy_J)f_J(x+c_1ty)
\\\qquad\qquad\qquad\qquad\qquad\qquad\qquad+t y_Iy_J(g_J(x+c_1ty)+c_1 y_K \nabla_Kf_J(x+c_1ty))\Bigr]\:d\sigma(y)
\\
-\int_{c_2t}^{c_1t} \int_{S^2} r^{-1} (\delta_{IJ}-3y_Iy_J)(tg_J(x+ry)+f_J(x+ry))\:d\sigma(y)\:dr.
\end{multline}
After a simple change of variables, this is
\begin{multline}
\label{fund.soln.2}
4\pi v_I(t,x)=\int_{S^2} \Bigl[(2\delta_{IJ} - 4 y_I y_J)f_J(x+c_2ty)
\\\qquad\qquad\qquad\qquad\qquad\qquad\qquad
+t (\delta_{IJ}-y_Iy_J)(g_J(x+c_2ty)+c_2 y_K \nabla_Kf_J(x+c_2ty))\Bigr]\:d\sigma(y)
\\\!\!\!\!\!\!\!\!\!\!\!\!\!\!\!\!+\int_{S^2} \Bigl[(-\delta_{IJ}+4y_Iy_J)f_J(x+c_1ty)
\\\qquad\qquad\qquad\qquad\qquad\qquad\qquad+t y_Iy_J(g_J(x+c_1ty)+c_1 y_K \nabla_Kf_J(x+c_1ty))\Bigr]\:d\sigma(y)
\\
-\int_{c_2}^{c_1} \int_{S^2} c^{-1} (\delta_{IJ}-3y_Iy_J)(tg_J(x+cty)+f_J(x+cty))\:d\sigma(y)\:dc.
\end{multline}
Using this representation, we shall be able to prove estimates on $v$ using techniques similar to those
for wave equations.

\subsection{Pointwise estimates in $\R^3$}
In \cite{KSS3}, the authors adapted a $L^1-L^\infty$ H\"ormander-type estimate to eliminate the dependence
on Lorentz invariance.  As such, these estimates could be applied to study multiple speed wave equations
and certain boundary value problems.  In this section, we prove analogues of these pointwise estimates.

The first of these estimates is for the homogeneous equation \eqref{boundaryless.equation}.
\begin{lemma}
\label{lemma.ptwise.bdyless.data}
Let $v$ be a solution to \eqref{boundaryless.equation}.  Then,
\begin{equation}
\label{ptwise.bdyless.data}
(1+t+|x|)|v(t,x)|\lesssim \sum_{|\alpha|\le 4} \|\la x\ra^{|\alpha|} \partial^\alpha f\|_2
+\sum_{|\alpha|\le 3} \|\la x\ra^{1+|\alpha|}\partial^\alpha g\|_2.
\end{equation}
\end{lemma}

\noindent{\em Proof of Lemma \ref{lemma.ptwise.bdyless.data}:}
Using the following estimates for the wave equation from \cite{MNS2}\footnote{See the proof of
Lemma 2.2.}
\begin{equation}
\label{wave.data.1}
(1+t+|x|) t\int_{S^2} |h(x+ty)|\:d\sigma(y) \lesssim \sum_\ss{|\alpha|+\mu\le 3\\\mu\le 1}
\int_{\R^3} |(|z|\partial_{|z|})^\mu Z^\alpha h(z)|\:\frac{dz}{\la z\ra}
\end{equation}
and
\begin{equation}
\label{wave.data.2}
(1+t+|x|)\int_{S^2} |h(x+ty)|\:d\sigma(y)\lesssim \sum_\ss{|\alpha|+\mu\le 3\\\mu\le 1}
\int_{\R^3} |(|z|\partial_{|z|})^\mu Z^\alpha h(z)|\:\frac{dz}{\la z\ra^2},
\end{equation}
\eqref{ptwise.bdyless.data} follows from \eqref{fund.soln.2} and the Schwarz inequality.\qed

We now prove the corresponding result for the inhomogeneous equation.
\begin{lemma}
\label{lemma.ptwise.bdyless.inhom}
Let $w$ solve $L w=G$ for $(t,x)\in \R_+\times\R^3$, and suppose that
$w(t,x)=0$ for $t\le 0$.  Then,
\begin{equation}
\label{ptwise.bdyless.inhom}
(1+t+|x|)|w(t,x)|\lesssim  \sum_\ss{|\alpha|+\mu\le 3\\\mu\le 1} \int_0^t \int_{\R^3} |S^\mu Z^\alpha
G(s,y)|\:\frac{dy\:ds}{|y|}.
\end{equation}
Suppose further that $G(s,y)=0$ when $|y|>20c_1 s$.  Then,
\begin{equation}
\label{ptwise.bdyless.inhom.localized}
(1+t)|w(t,x)|\lesssim \sum_\ss{|\alpha|+\mu\le 3\\\mu\le 1} \int_{\theta t}^t \int_{\R^3} |S^\mu Z^\alpha
G(s,y)|\:\frac{dy\:ds}{|y|},\quad \text{for } |x|<c_2 t/10,
\end{equation}
where $\theta$ is some constant depending on $c_2$.
\end{lemma}

\noindent{\em Proof of Lemma \ref{lemma.ptwise.bdyless.inhom}:}  This follows as in the proof of the
previous lemma, using the following bound from \cite{KSS3}\footnote{See the proof of Proposition 2.1.}
$$\int_0^t \int_{S^2} (t-s) |F(s,x-(t-s)y)|\:d\sigma(y)\:ds
\lesssim  \sum_\ss{|\alpha|+\mu\le 3\\\mu\le 1} \int_0^t \int_{\R^3} |S^\mu Z^\alpha F(s,y)|\:
\frac{dy\:ds}{|y|}.$$
The lemma then follows from \eqref{fund.soln.2} and Duhamel's principle.

\begin{figure}[h]
\begin{center}
\epsfxsize=15cm \epsffile{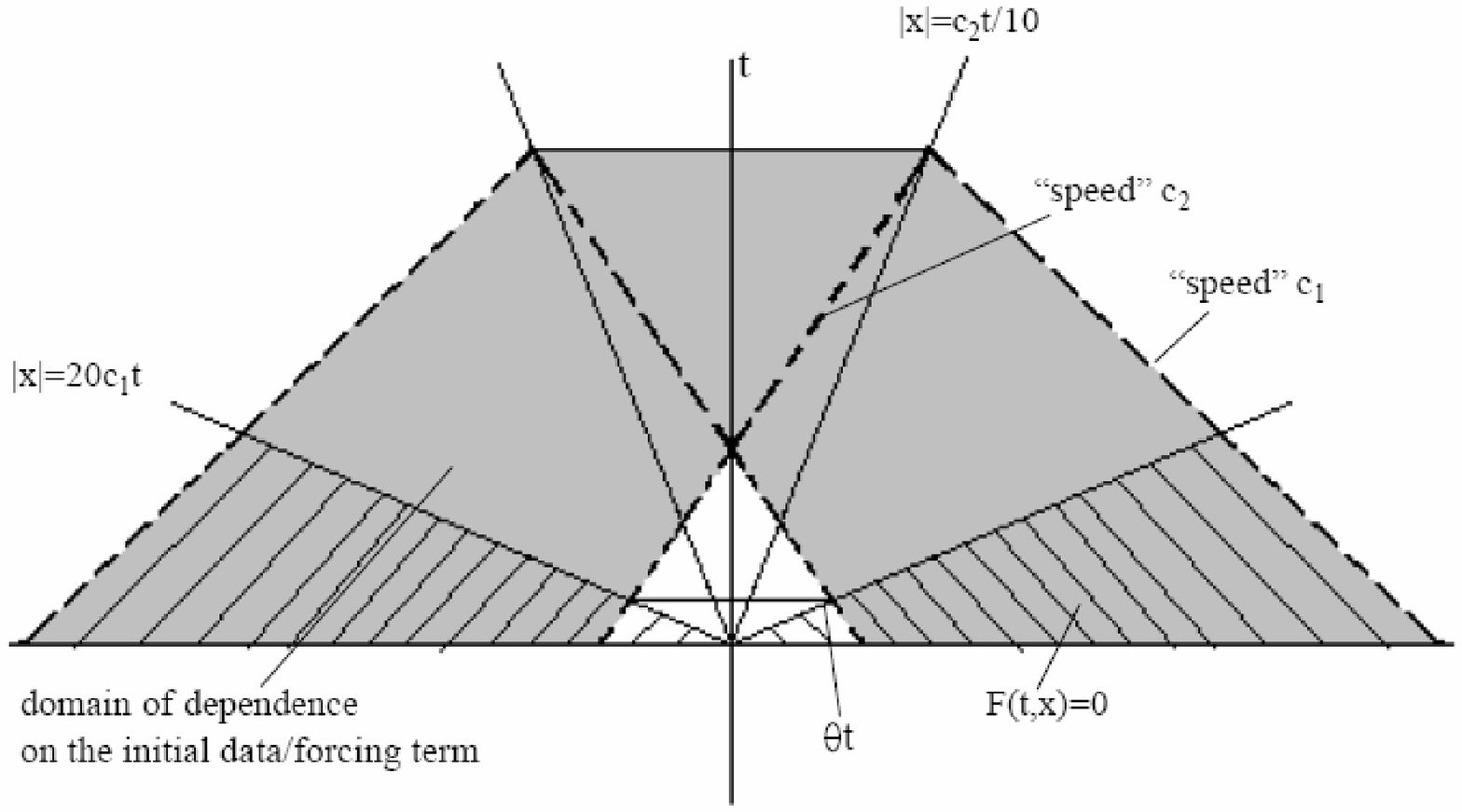}
\end{center}
\caption{An illustration of the definition of $\theta t$.} \label{trsgsvtime}
\end{figure}

The second estimate follows from the appropriate Huygens' principle for $L$.  See Fig. \ref{trsgsvtime}.\qed

By Duhamel's principle and the positivity of \eqref{fund.soln.2}, we also have
\begin{lemma}
\label{lemma.ptwise.bdyless.2}\footnote{We note that a portion
of this
lemma has previously appeared in \cite{XQ}.}
Let $v$ solve \eqref{boundaryless.equation}.  Then,
\begin{multline}
\label{ptwise.bdyless.data.2}
|x||v(t,x)|\lesssim \sum_{i=1,2} \int_{||x|-c_i(t-s)|}^{|x|+c_i(t-s)}
\sup_{|\theta|=1} |f(\rho\theta)|\:d\rho
+\int_{c_2}^{c_1} \int_{||x|-c(t-s)|}^{|x|+c(t-s)} \sup_{|\theta|=1} |f(\rho\theta)|\:d\rho\:dc
\\+\sum_{i=1,2} \int_{||x|-c_i(t-s)|}^{|x|+c_i(t-s)} \sup_{|\theta|=1} |\nabla f(\rho\theta)|
\rho\:d\rho
\\+\sum_{i=1,2} \int_{||x|-c_i(t-s)|}^{|x|+c_i(t-s)} \sup_{|\theta|=1} |g(\rho\theta)|\rho\:d\rho
+\int_{c_2}^{c_1} \int_{||x|-c(t-s)|}^{|x|+c(t-s)} \sup_{|\theta|=1} |g(\rho\theta)|\rho\:d\rho\:dc.
\end{multline}
Moreover, let $w$ solve $L w = G$ for $(t,x)\in \R_+\times\R^3$, and suppose that $w(t,x)=0$ for $t\le 0$.  Then
\begin{multline}
\label{ptwise.bdyless.inhom.2}
|x||w(t,x)|\lesssim \sum_{i=1,2}\int_0^t \int_{||x|-c_i(t-s)|}^{|x|+c_i(t-s)} 
\sup_{|\theta|=1} |G(s,\rho\theta)|\rho
\:d\rho\:ds
\\+\int_{c_1}^{c_2} \int_0^t \int_{||x|-c(t-s)|}^{|x|+c(t-s)} \sup_{|\theta|=1} |G(s,\rho\theta)|\rho\:
d\rho\:ds\:dc.
\end{multline}
\end{lemma}

\subsection{Pointwise estimates in $\ext$}
Here, we wish to establish the main decay estimates which shall be used in the sequel.  
At the first order of difficulty, these are exterior domain analogs of \eqref{ptwise.bdyless.data}
and \eqref{ptwise.bdyless.inhom}.  For these estimates, we shall look at solutions to the
following
\begin{equation}
\label{ptwise.equation}
\begin{cases}
Lu(t,x)=F(t,x),\quad (t,x)\in \R_+\times\ext,\\
u(t,x)|_{\bdy}=0,\\
u(0,x)=f(x),\quad \partial_t u(0,x)=g(x),\\
\supp f,g \subset \{|x|\ge 6\}.
\end{cases}
\end{equation}
For such solutions, we shall have the following estimates.
\begin{theorem}
\label{theorem.ptwise}
Let $\K\subset\R^3$ be a bounded, nontrapping obstacle with smooth boundary.  Then smooth
solutions of \eqref{ptwise.equation} satisfy
\begin{multline}
\label{ptwise}
(1+t+|x|)|S^\nu Z^\alpha u(t,x)|
\lesssim \sum_\ss{j+|\beta|+k\le M+\nu+7\\j\le 1} \|\la x\ra^{j+|\beta|} \nabla^\beta
\partial_t^{k+j} u(0,\cd)\|_2
\\+ \int_0^t \int_\ext \sum_\ss{|\beta|+\mu\le M+\nu+6\\\mu\le\nu+1} |S^\mu Z^\beta F(s,y)|\:\frac{dy\:ds}{|y|}
\\+\int_0^t \sum_\ss{|\beta|+\mu\le M+\nu+3\\\mu\le\nu+1} \|S^\mu \partial^\beta F(s,\cd)\|_{L^2(\{|x|<4\})}\:ds
\end{multline}
for any $|\alpha|=M$ and any $\nu$.
\end{theorem}

\noindent{\em Proof of Theorem \ref{theorem.ptwise}:}  This proof resembles those in \cite{KSS3} and \cite{MS1}
for the wave equation quite closely.  A portion of this argument was also given previously in \cite{XQ}.

As a first reduction, we prove that
\begin{multline}
\label{ptwise.1}
(1+t+|x|)|S^\nu Z^\alpha u(t,x)|\lesssim \sum_\ss{j+|\beta|+k\le M+\nu+4\\j\le 1}
\|\la x\ra^{j+|\beta|}\nabla^\beta \partial_t^{k+j} u(0,\cd)\|_2
\\+\int_0^t \int_\ext \sum_\ss{|\beta|+\mu\le |\alpha|+\nu+3\\\mu\le\nu+1} |S^\mu Z^\beta F(s,y)|
\:\frac{dy\:ds}{|y|}
\\+\sup_{|y|\le 2, 0\le s\le t} (1+s) \sum_\ss{|\beta|+\mu\le |\alpha|+\nu\\\mu\le\nu}
\Bigl(|S^\mu Z^\beta u'(s,y)| + |S^\mu Z^\beta u(s,y)|\Bigr). 
\end{multline}

\noindent{\em Proof of \eqref{ptwise.1}:}
To do this, we follow the proof of Lemma 4.2 of \cite{KSS3}.  The estimate is obvious
for $|x|<2$.  Thus, for the remainder of the proof, we shall assume that $|x|\ge 2$.  As such,
for a smooth $\rho$ satisfying $\rho(r)\equiv 1$ for $r\ge 2$ and $\rho(r)\equiv 0$ for $r\le 1$, it
suffices to establish the estimate for $w(t,x)=\rho(|x|) S^\nu Z^\alpha u(t,x)$, which solves
the boundaryless\footnote{Recall that we are assuming $\K\subset\{|x|<1\}$.} equation
\begin{multline*}
Lw=\rho L S^\nu Z^\alpha u - 2c_2^2\nabla \rho \cdot \nabla (S^\nu Z^\alpha u) - c_2^2(\Delta \rho) S^\nu Z^\alpha u
\\-(c_1^2-c_2^2)\nabla(\nabla\rho\cdot S^\nu Z^\alpha u) - (c_1^2-c_2^2)\nabla\rho \nabla
\cdot (S^\nu Z^\alpha u).
\end{multline*}
We write $w=w_1+w_2$ where $Lw_1=\rho L S^\nu Z^\alpha u$ with $(w_1,\partial_t w_1)(0,\cd)=
(f,g)$.  The estimate for $w_1$ follows from \eqref{ptwise.bdyless.data} and \eqref{ptwise.bdyless.inhom}.

By \eqref{ptwise.bdyless.inhom.2}, we see that 
\begin{multline}\label{ptwise.w0}
|w_2(t,x)|\lesssim \frac{1}{|x|}\sum_{i=1,2}\int_0^t \int_{|r-c_i(t-s)|}^{r+c_i(t-s)} \sup_{|\theta|=1}
|G(s,\rho\theta)|\rho\:d\rho\:ds
\\+\frac{1}{|x|}\int_{c_2}^{c_1} \int_0^t \int_{|r-c(t-s)|}^{r+c(t-s)} \sup_{|\theta|=1}
|G(s,\rho\theta)|\rho\:d\rho\:ds\:dc
\end{multline}
where
$$
G(t,x)=Lw-\rho LS^\nu Z^\alpha u.$$  
We shall now work with the last term in \eqref{ptwise.w0}.  The bound for the other term follows similarly.
By noting that $G$ vanishes for $|x|\ge 2$, it follows that we must have
$$-2\le |x|-c(t-s)\le 2$$
for the $d\rho$ integral to be nonzero.  That is, we must have
$$\frac{ct-|x|-2}{c}\le s\le \frac{ct-|x|+2}{c}.$$
Thus, we conclude that the last term in \eqref{ptwise.w0} is
\begin{multline*}
\lesssim \frac{1}{|x|}\int_{c_2}^{c_1} \frac{1}{1+|ct-|x||}\\\times\sum_\ss{|\beta|+\mu\le M+\nu\\\mu\le\nu}
\sup_\ss{\frac{ct-|x|-2}{c}\le s\le \frac{ct-|x|+2}{c}\\
|y|\le 2} (1+s)\Bigl(|S^\mu Z^\beta u'(s,y)| + |S^\mu Z^\beta u(s,y)|\Bigr)\:dc
\end{multline*}
from which the bound for $|w_2|$ follows.\qed

In order to complete the proof of Theorem \ref{theorem.ptwise}, it will suffice to show that
$$(1+t)\sum_\ss{|\beta|+\mu\le M+\nu+1\\\mu\le\nu}\sup_{|x|<2} |S^\mu \partial^\beta u(t,x)|$$
is bounded by the right side of \eqref{ptwise}.  Using smooth cutoffs, we may examine the following
cases separately:
\begin{enumerate}
\item[{\em Case 1:}] $u$ has vanishing data and $F(s,y)$ vanishes for $|y|>4$
\item[{\em Case 2:}] $F(s,y)=0$ for $|y|<3$.
\end{enumerate}

We use the following consequence of the Fundamental Theorem of Calculus and Sobolev embedding
\begin{equation}
\label{FTC}
\begin{split}
(1+t)\sum_\ss{|\beta|+\mu\le M+\nu+1\\\mu\le\nu} \sup_{|x|<2} &|S^\mu \partial^\beta u(t,x)|\\
&\lesssim \int_0^t \sum_\ss{|\beta|+\mu\le M+\nu+2\\\mu\le \nu, j\le 1} \|(s\partial_s)^j 
S^\mu \partial^\beta u'(s,\cd)\|_{L^2(\{|x|<4\})}\:ds\\
&\lesssim \int_0^t \sum_\ss{|\beta|+\mu\le M+\nu+3\\\mu\le\nu+1}\|S^\mu \partial^\beta
u'(s,\cd)\|_{L^2(\{|x|<4\})}\:ds.
\end{split}
\end{equation}
Here, for the first inequality, we have also used the fact that the Dirichlet boundary conditions
permit us to bound $u$ locally by $u'$.
The bound in {\em Case 1} now follows from an application of \eqref{local.energy.high}.

To complete the proof, we need only examine {\em Case 2}.  Here, we write $u=w+u_r$ where
$w$ solves the boundaryless equation $Lw=F$ with $(w,\partial_t w)(0,\cd)=(u,\partial_t u)(0,\cd)$.
We fix $\eta\in C^\infty(\R^3)$ with $\eta(x)\equiv 1$ for $|x|<2$ and $\eta(x)\equiv 0$ for $|x|\ge 3$.
Letting $\tilde{u}=\eta w + u_r$, we notice that $u=\tilde{u}$ for $|x|<2$, that $\tilde{u}$ solves
\begin{equation}\label{tildeu.equation}
L\tilde{u} = -2c_2^2\nabla \eta\cdot\nabla w - c_2^2(\Delta\eta) w - (c_1^2-c_2^2)
\nabla(\nabla\eta\cdot w) - (c_1^2-c_2^2)\nabla\eta \nabla\cdot w,
\end{equation}
and that $\tilde{u}$ has vanishing Cauchy data.

As the right side vanishes unless $2\le |x|\le 3$, we may apply the result of the previous case to see that
$$
(1+t)\sum_\ss{|\beta|+\mu\le M+\nu+1\\\mu\le\nu} \sup_{|x|<2} |S^\mu \partial^\beta \tilde{u}(t,x)|
\lesssim \int_0^t \sum_\ss{|\beta|+\mu\le M+\nu+4\\\mu\le\nu} \|S^\mu \partial^\beta w(s,\cd)\|_{L^\infty(2\le
|x|\le 3)}\:ds.
$$

It remains only to bound the term on the right of the preceding equation.  To do this, we apply
\eqref{ptwise.bdyless.data.2}, \eqref{ptwise.bdyless.inhom.2}, and Sobolev's lemma on
$S^2$.  For the forcing term, for example, 
we have
$$\sum_\ss{|\beta|+\mu\le M+\nu+6\\\mu\le \nu}\int_0^t \int_{c_2}^{c_1} \int_0^s \int_{|c(s-\tau)-|y||\le 4}
|S^\mu Z^\beta F(\tau,y)|\:\frac{dy\:d\tau}{|y|}\:dc\:ds$$
corresponding to the last term in \eqref{ptwise.bdyless.inhom.2}.  For fixed $c$, we note that
the sets $\Lambda_s = \{(\tau,y)\,:\, 0\le \tau\le s,\,|c(s-\tau)-|y||\le 4\}$ have finite overlap.  I.e.,
for $|s-s'|\ge 10$, $\Lambda_s\cap\Lambda_{s'}=\emptyset$.  Thus, we see that the desired bound holds.  
Similar arguments can be used to estimate the remaining terms
in \eqref{ptwise.bdyless.data.2} and \eqref{ptwise.bdyless.inhom.2}, which completes the proof.
\qed

We will also require the following pointwise estimate which follows from arguments similar to those in
\cite{MNS1}.
\begin{theorem}
\label{theorem.ptwise.2}
Suppose $\K\subset \{|x|<1\}\subset \R^3$ is a nontrapping obstacle with smooth boundary.
Let $u$ solve
\begin{equation}
\label{no.data.equation}
\begin{cases}
Lu(t,x)=F(t,x),\quad (t,x)\in \R_+\times\ext,\\
u|_\bdy=0,\\
u(t,x)=0,\quad t\le 0.
\end{cases}
\end{equation}
Suppose further that $F(t,x)=0$ when $|x|>20c_1t$.  Then, if $|x|<c_2t/10$ and $t>1$,
\begin{multline}
\label{ptwise.2}
(1+t+|x|)|S^\nu Z^\alpha u(t,x)| \lesssim \sum_\ss{|\beta|+\mu\le |\alpha|+\nu+6\\\mu\le\nu+1}
\int_{\theta t}^t \int_\ext |S^\mu Z^\beta F(s,y)|\:\frac{dy\:ds}{|y|}
\\+\sup_{0\le s\le t} (1+s) \sum_\ss{|\beta|+\mu\le |\alpha|+\nu+2\\\mu\le\nu} \|S^\mu \partial^\beta
F(s,\cd)\|_2.
\end{multline}
\end{theorem}

\noindent{\em Proof of Theorem \ref{theorem.ptwise.2}:}
By arguing as in \eqref{ptwise.1} using \eqref{ptwise.bdyless.inhom.localized}, we see that 
for $|x|<c_2 t/10$,
\begin{multline}
\label{ptwise.2.1}
(1+t)|S^\nu Z^\alpha u(t,x)|\lesssim \int_{\theta t}^t \int_\ext 
\sum_\ss{|\beta|+\mu\le |\alpha|+\nu+3\\\mu\le\nu+1} |S^\mu Z^\beta F(s,y)|\:\frac{dy\:ds}{|y|}
\\+\sup_{|y|\le 2, 0\le s\le t} (1+s) \sum_\ss{|\beta|+\mu\le |\alpha|+\nu+1\\\mu\le\nu}
|S^\mu \partial^\beta u(s,y)|.
\end{multline}
Thus, we need only show that the last term is controlled by the right side of \eqref{ptwise.2}.

First suppose that $F(s,y)=0$ if $|y|>4$.  Then, by a Sobolev estimate and \eqref{local.energy.high},
it follows that
\begin{equation}
\label{ptwise.2.2}
(1+t)\sup_{|y|<2} \sum_\ss{|\beta|+\mu\le |\alpha|+\nu+1\\\mu\le\nu}|S^\mu \partial^\beta u(t,x)|
\lesssim \sup_{0\le s\le t} \sum_\ss{|\beta|+\mu\le |\alpha|+\nu+2\\\mu\le\nu} (1+s)\|S^\mu \partial^\beta
F(s,\cd)\|_2.
\end{equation}
Here we are also using that the Dirichlet boundary conditions allow us to control $u$ locally 
by $u'$.  Thus, for the remainder of the proof, we may assume that $F(s,y)=0$ for $|y|\le 3$.

Letting $\tilde{u}$ be as in the proof of Theorem \ref{theorem.ptwise}, we see from \eqref{tildeu.equation} and
\eqref{ptwise.2.2} that
\begin{multline}
\label{ptwise.2.3}
(1+t)\sup_{|y|<2} \sum_\ss{|\beta|+\mu\le |\alpha|+\nu+1\\\mu\le\nu} |S^\mu \partial^\beta u(t,x)|
\\\lesssim \sup_{0\le s\le t} \sum_\ss{|\beta|+\mu\le |\alpha|+\nu+3\\\mu\le\nu} (1+s)\|S^\mu 
\partial^\beta w(s,\cd)\|_{L^\infty(2\le |x|\le 3)},
\end{multline}
where $w$ is a solution to the boundaryless equation $Lw=F$ with vanishing initial data.
Provided that $s>30/c_2$, the theorem follows from \eqref{ptwise.bdyless.inhom.localized},
\eqref{ptwise.2.1}, \eqref{ptwise.2.2}, and
\eqref{ptwise.2.3}.  Else, when $s\le 30/c_2$, we may use
\eqref{ptwise.bdyless.inhom}
and finite propagation speed to bound this term by the second term in the right of \eqref{ptwise.2}, which
completes the proof.\qed

\subsection{Boundary term estimates}\label{boundary}
In order to handle the boundary terms that arise in \eqref{energy.L.d}, we shall give an estimate which is in the 
spirit of 
the original ones of \cite{MS1, MS2} for the wave equation.  This is the key estimate, along
with the decay of local energy \eqref{local.energy}, that allows us to drop the star-shapedness assumption on $\K$.
Here, we are merely isolating the proof of the bound for the right side of \eqref{FTC}.

\begin{lemma}
\label{lemma.bdy.term}  Let $\K$ be a bounded, nontrapping obstacle with smooth boundary.
Suppose that $u$ is a smooth solution of \eqref{no.data.equation}.
Then, 
\begin{multline}
\label{bdy.term}
\int_0^t \|S^\nu \partial^\alpha u'(s,\cd)\|_{L^2(|x|<1)}\:ds \lesssim \sum_\ss{|\beta|+\mu\le M+\nu+3\\\mu\le \nu}
\int_0^t \int_\ext |S^\mu Z^\beta F(s,y)|\:\frac{dy\:ds}{|y|}
\\+ \int_0^t \sum_\ss{|\beta|+\mu\le M+\nu\\\mu\le \nu} \|S^\mu \partial^\beta F(s,\cd)\|_{L^2(\{|x|<4\})}\:ds
\end{multline}
for any $|\alpha|=M$ and $\nu$.
\end{lemma}

\bigskip
\newsection{Sobolev-type estimates and null form estimates}

\subsection{Weighted Sobolev estimates}
In the sequel, we shall require the following rather standard weighted Sobolev estimate from 
\cite{K}.  This follows by applying Sobolev embedding on $\R\times S^2$.  The decay results
from the difference in the volume elements between $\R\times S^2$ and, say, $\R^3$.
\begin{lemma}
\label{lemma.weighted.Sobolev}
Suppose that $h\in C^\infty(\ext)$.  Then, for $R\ge 1$,
\begin{equation}
\label{weighted.Sobolev}
\|h\|_{L^\infty(R/2<|x|<R)} 
\lesssim R^{-1} \sum_{|\alpha|+|\beta|\le 2} \|\tilde{\Omega}^\alpha \partial_x^\beta h\|_{L^2(R/4<|x|<2R)}.
\end{equation}
\end{lemma}

\subsection{Klainerman-Sideris decay estimates}
In this section, we study another class of weighted Sobolev-type estimates.  These follow from the boundaryless
estimates of \cite{Si}, which are closely related to estimates that originally appeared in \cite{KS}.
The interested reader should also consult \cite{ST} for a unified approach to proving such estimates.

The main boundaryless estimate of \cite{Si}\footnote{See (3.24).} to consider states 
\begin{multline}
\label{KS.bdyless}
\|\la c_1t-r\ra P_1 \partial\nabla h(t,\cd)\|_2 + \|\la c_2t-r\ra P_2 \partial\nabla h(t,\cd)\|_2
\\\lesssim \sum_{|\alpha|+\mu\le 1} \|S^\mu Z^\alpha h'(t,\cd)\|_2
+ t\|Lh(t,\cd)\|_2.
\end{multline}
In order to establish similar bounds when there is a boundary, we use ideas from \cite{MNS1}.  This yields
\begin{lemma}
\label{lemma.KS.1}
Let $\K\subset\{|x|<1\}\subset\R^3$ be an obstacle with smooth boundary.
Suppose $u(t,x)\in C^\infty(\R_+\times\ext)$ and $u|_\bdy=0$.  Then,
\begin{multline}
\label{KS.1}
\|\la c_1t-r\ra P_1 \partial \nabla S^\nu Z^\alpha u(t,\cd)\|_2 + \|\la c_2t-r\ra P_2 \partial\nabla 
S^\nu Z^\alpha u(t,\cd)\|_2
\\\lesssim \sum_\ss{|\beta|+\mu\le M+\nu+1\\\mu\le\nu+1} \|S^\mu Z^\alpha u'(t,\cd)\|_2
+ t\sum_\ss{|\beta|+\mu\le M+\nu\\\mu\le\nu}\|S^\mu Z^\beta Lu(t,\cd)\|_2 \\+
(1+t)\sum_{\mu\le\nu} \|S^\mu u'(t,\cd)\|_{L^2(|x|<2)} 
\end{multline}
for any fixed $|\alpha|=M$ and $\nu$.
\end{lemma}

\noindent{\em Proof of Lemma \ref{lemma.KS.1}:} We first show that
\begin{multline}
\label{KS.1.1}
\sum_{\kappa=1,2} \|\la c_\kappa t-r\ra P_\kappa \partial \nabla S^\nu Z^\alpha u(t,\cd)\|_2
\lesssim \sum_\ss{|\beta|+\mu\le M+\nu+1\\\mu\le\nu+1} \|S^\mu Z^\beta u'(t,\cd)\|_2
\\+ \sum_\ss{|\beta|+\mu\le M+\nu\\\mu\le\nu} t \|S^\mu Z^\beta Lu(t,\cd)\|_2
+ (1+t)\sum_\ss{|\beta|+\mu\le M+\nu+2\\\mu\le\nu} \|S^\mu \partial^\beta u(t,\cd)\|_{L^2(|x|<3/2)}.
\end{multline}
This estimate is trivial if the norm in the left is over $\{|x|<3/2\}$ as the coefficients of $Z$ are $O(1)$ on
this set.  To handle the case when the norms are over $\{|x|\ge 3/2\}$, we fix $\eta\in C^\infty(\R^3)$ with
$\eta(x)\equiv 0$ for $|x|<1$ and $\eta(x)\equiv 1$ for $|x|>3/2$ and apply \eqref{KS.bdyless} to 
$h(t,x)=\eta(x)u(t,x)$ which solves the boundaryless equation
$$Lh=\eta(x)Lu - 2c_2^2 \nabla \eta\cdot\nabla u - c_2^2(\Delta \eta)u - (c_1^2-c_2^2)\nabla(\nabla\eta\cdot
u)-(c_1^2-c_2^2)\nabla\eta\nabla\cdot u.$$
This proves \eqref{KS.1.1} as the last four terms are supported in $\{|x|\le 3/2\}$.

It remains to show that the last term in \eqref{KS.1.1} is bounded by the right side of \eqref{KS.1}.  Here
we simply apply (a trivial modification of) \eqref{ell.reg.local} and the differential inequality
$$\sum_{\kappa=1,2}|c_\kappa t-r||P_\kappa\partial_t \nabla v(t,x)|
\lesssim \sum_{|\alpha|+\mu\le 1} |S^\mu Z^\alpha \nabla v(t,x)|+t|Lv(t,x)|,$$
which is also from \cite{Si}.\footnote{See (3.10b).} \qed

We also have the following variants of the above.
\begin{corr}
\label{corr.ks.1}
Let $\K\subset\{|x|<1\}\subset\R^3$ be an obstacle with smooth boundary.
Suppose $u(t,x)\in C^\infty(\R_+\times\ext)$ and $u|_\bdy=0$.  Then,
\begin{multline}
\label{KS.2}
|x|\la c_\kappa t-r\ra |P_\kappa \partial\nabla S^\nu Z^\alpha u(t,x)|
\lesssim \sum_\ss{|\beta|+\mu\le M+\nu+3\\\mu\le \nu+1} \|S^\mu Z^\beta u'(t,\cd)\|_2
\\+ \sum_\ss{|\beta|+\mu\le M+\nu+2\\\mu\le \nu} t\|S^\mu Z^\beta Lu(t,\cd)\|_2
+ (1+t)\sum_{\mu\le\nu} \|S^\mu u'(t,\cd)\|_{L^2(|x|<2)},
\end{multline}
\begin{multline}
\label{KS.3}
\|\la c_\kappa t-r\ra P_\kappa \partial\nabla S^\nu Z^\alpha u(t,\cd)\|_{L^2(|x|>t/4)}
\lesssim \sum_\ss{|\beta|+\mu\le M+\nu+1\\\mu\le\nu+1} \|S^\mu Z^\beta u'(t,\cd)\|_{L^2(|x|>t/8)}
\\+ (1+t)\sum_\ss{|\beta|+\mu\le M+\nu\\\mu\le\nu} \|S^\mu Z^\beta Lu(t,\cd)\|_{L^2(|x|>t/8)},
\end{multline}
and
\begin{multline}
\label{KS.4}
\sup_{|x|>t/2} |x|\la c_\kappa t-r\ra |P_\kappa \partial\nabla S^\nu Z^\alpha u(t,\cd)|
\lesssim \sum_\ss{|\beta|+\mu\le |\alpha|+\nu+3\\\mu\le\nu+1} \|S^\mu Z^\beta u'(t,\cd)\|_{L^2(|x|>t/8)}
\\+(1+t)\sum_\ss{|\beta|+\mu\le |\alpha|+\nu+2\\\mu\le\nu} \|S^\mu Z^\beta Lu(t,\cd)\|_{L^2(|x|>t/8)}
\end{multline}
for any $\kappa=1,2$ and fixed $|\alpha|=M$ and $\nu$.
\end{corr}

\noindent{\em Proof of Corollary \ref{corr.ks.1}:}
Estimates \eqref{KS.2} and \eqref{KS.4} follow from \eqref{KS.1} and \eqref{KS.3} respectively using
\eqref{weighted.Sobolev}.  Estimate \eqref{KS.3} follows from applying \eqref{KS.bdyless} to
$h(t,x)=\eta(x/\la t\ra) u(t,x)$ where $\eta(z)\equiv 1$ for $|z|>1/4$ and $\eta(z)\equiv 0$ for $|z|<1/8$.\qed

We also have
\begin{lemma}
\label{lemma.KS.5}
Let $\K\subset\{|x|<1\}\subset\R^3$ be an obstacle with smooth boundary.  Suppose $u(t,x)\in C^\infty(\R_+\times
\ext)$ and $u|_\bdy=0$.  Then,
\begin{multline}
\label{KS.5}
\la r\ra \la c_\kappa t-r\ra^{1/2} |P_\kappa \partial S^\nu Z^\alpha u(t,x)|
\lesssim \sum_\ss{|\beta|+\mu\le M+\nu+2\\\mu\le\nu+1} \|S^\mu Z^\beta u'(t,\cd)\|_2
\\+(1+t)\sum_\ss{|\beta|+\mu\le M+\nu+1\\\mu\le\nu} \|S^\mu Z^\beta Lu(t,\cd)\|_2
+(1+t)\sum_{\mu\le\nu} \|S^\mu u'(t,\cd)\|_{L^2(|x|<2)}
\end{multline}
and
\begin{multline}
\label{KS.6}
\sup_{|x|>t/4} \la r\ra\la c_\kappa t-r\ra^{1/2} |P_\kappa \partial S^\nu Z^\alpha u(t,x)|
\lesssim \sum_\ss{|\beta|+\mu\le M+\nu+2\\\mu\le\nu+1} \|S^\mu Z^\beta u'(t,\cd)\|_{L^2(|x|>t/8)}
\\+(1+t)\sum_\ss{|\beta|+\mu\le M+\nu+1\\\mu\le\nu} \|S^\mu Z^\beta Lu(t,\cd)\|_{L^2(|x|>t/8)}
\end{multline}
for any $\kappa=1,2$ and fixed $|\alpha|=M$ and $\nu$.
\end{lemma}

\noindent{\em Proof of Lemma \ref{lemma.KS.5}:}
Using the arguments of the proof of Lemma \ref{lemma.KS.1}, this follows from the boundaryless 
estimate
$$
\la r\ra\la c_\kappa t-r\ra^{1/2} |P_\kappa \partial h(t,x)|
\lesssim \sum_\ss{|\beta|+\mu\le 2\\\mu\le 1} \|S^\mu Z^\beta h'(t,\cd)\|_2
+ \la t\ra \sum_\ss{|\beta|\le 1} \|Z^\beta Lh(t,\cd)\|_2 
$$
of \cite{Si}.\footnote{See (3.20c) and (3.24).}\qed

\subsection{Null form estimates}
The additional decay afforded to us by the null condition is encapsulated in the following
estimate of \cite{Si}.\footnote{See Proposition 3.2.}  
This is closely related to those for multiple speed systems of wave
equations used, e.g., in \cite{Si3}.
\begin{lemma}
\label{lemma.null.condition}
Assume that $Q$ satisfies the null condition \eqref{null.condition}.  Let 
$\mathcal{N}=\{(\alpha,\beta,\gamma)\neq (1,1,1), (2,2,2)\}$ be the set of
nonresonant indices.  Then,
\begin{multline}
\label{null.condition.decay}
|\la u, Q(v,w)\ra|\\\lesssim \frac{1}{r} |u|\sum_{|a|\le 1} \Bigl[|\nabla \tilde{\Omega}^a v|
|\nabla w| + |\nabla\tilde{\Omega}^a w||\nabla v|
+ |\nabla^2 v||\tilde{\Omega}^a w| + |\nabla^2 w||\tilde{\Omega}^a v|\Bigr]
\\+\sum_{\mathcal{N}} |P_\alpha u|\Bigl[
|P_\beta \nabla^2 v| |P_\gamma \nabla w| + |P_\beta \nabla^2 w||P_\gamma\nabla v|\Bigr],
\end{multline}
for sufficiently regular $u,v,w$.
\end{lemma}

\bigskip
\newsection{Proof of global existence}

In this section, we prove Theorem \ref{main.theorem}.  We take $N=112$ in \eqref{data.smallness}, but this is far from
optimal.  By scaling in the $t$ variable, we may take $c_1=1$ without loss of generality.

We begin by making a reduction that allows us to avoid technicalities related to the
compatibility conditions.  While the reduction from \cite{KSS3}\footnote{See also \cite{MS1} and 
\cite{MNS1}.} works when we have truncated at the quadratic level, the necessary scaling breaks down
when the higher order terms are present.  We instead use a reduction that is more reminiscent of
that from \cite{MNS2}.

We begin by noting that if $\varepsilon$ in \eqref{data.smallness} is sufficiently small, then there is
a constant $C_0$ for which
\begin{equation}
\label{local.solution.estimate}
\sup_{0\le t\le 2} \sum_{|\alpha|\le 112} \|\partial^\alpha u(t,\cd)\|_{L^2(|x|\le 25)}\le C_0\varepsilon.
\end{equation}
This follows from well-known local existence theory.  See, e.g., \cite{KSS1}, which is only stated for
diagonal wave equations but as the proofs only rely on energy estimates the results carry over to
the current setting.

On the other hand, over $\{|x|>5 (t+1)\}$, $u$ corresponds to a boundaryless solution.\footnote{Recall that
$\K\subset\{|x|<1\}$ by assumption.}  Thus by the estimates
which correspond to those that follow for the boundaryless problem\footnote{See, e.g., \cite{Si} and \cite{MNS2}.},
we have
\begin{multline}
\label{bdyless.solution.estimate}
\sup_{0\le t<\infty} \sum_{|\alpha|+\mu\le 111} \|S^\mu Z^\alpha u'(t,\cd)\|_{L^2(|x|\ge 5(t+1))}
\\+ \sup_\ss{|x|\ge 5(t+1)\\0\le t<\infty} (1+t+|x|)\sum_{|\alpha|+\mu\le 108} |S^\mu Z^\alpha
u'(t,x)|\le C_1\varepsilon.
\end{multline}

We fix a smooth cutoff function $\eta$ with $\eta(t,x)\equiv 1$ if $t\le 3/2$ and $|x|\le 20$,
$\eta(t,\cd)\equiv 0$ for $t>2$, and $\eta(\cd,x)\equiv 0$ for $|x|>25$.  Setting $u_0=\eta u$,
it follows that $u$ solves \eqref{main.equation} for $0<t<T$ if and only if $w=u-u_0$ solves
\begin{equation}
\label{w.equation}
\begin{cases}
Lw = (1-\eta)Q(\nabla u,\nabla^2 u) - [L,\eta]u, \quad (t,x)\in \R_+\times\ext,\\
w|_\bdy=0,\\
w(0,\cd)=(1-\eta)(0,\cd)f,\\
\partial_t w(0,\cd)=(1-\eta)(0,\cd)g - \eta_t(0,\cd)f
\end{cases}
\end{equation}
over the same interval.

We now fix another smooth cutoff $\beta$ with $\beta(z)\equiv 1$ for $z\ge 10$ and $\beta(z)\equiv 0$ for 
$z\le 6$.  Then, let $v$ solve the linear equation
\begin{equation}
\label{v.equation}
\begin{cases}
Lv = \beta\bigl(\frac{|x|}{t+1}\bigr)(1-\eta) Q(\nabla u, \nabla^2 u) - [L,\eta]u,\\
v|_\bdy = 0,\\
v(0,\cd)=(1-\eta)(0,\cd)f,\\
\partial_t v(0,\cd)=(1-\eta)(0,\cd)g - \eta_t(0,\cd)f.
\end{cases}
\end{equation}
We shall show that 
\begin{multline}
\label{v.estimate}
(1+t+|x|) \sum_{|\alpha|+\mu\le 102} |S^\mu Z^\alpha v(t,x)| + \sum_{|\alpha|+\mu\le 100} \|S^\mu Z^\alpha
v'(t,\cd)\|_2
\\+(\log(2+t))^{-1} \sum_{|\alpha|+\mu\le 98} \|\la x\ra^{-1/2} S^\mu Z^\alpha v'\|_{L^2_tL^2_x(S_t)}\le
C_2\varepsilon
\end{multline}
for some absolute constant $C_2>0$ and for any $t>0$.

\noindent{\em Proof of \eqref{v.estimate}:}
We first notice that by \eqref{ptwise} and \eqref{data.smallness} the first term in the left side of \eqref{v.estimate}
is
\begin{multline}
\label{v.equation.1}
\lesssim \varepsilon + \int_0^t \int_\ext \sum_{|\alpha|+\mu\le 108} |S^\mu Z^\alpha 
\beta\bigl(\frac{|x|}{t+1}\bigr)(1-\eta) Q(\nabla u,\nabla^2u)|\:\frac{dy\:ds}{|y|}
\\+ \int_0^t \int_\ext \sum_{|\alpha|+\mu\le 108} |S^\mu Z^\alpha [L,\eta]u|\:\frac{dy\:ds}{|y|}.
\end{multline}
Here, we have also used that Sobolev's lemma and the assumption that $0\in \K$ allow us to control
the last term in the right of \eqref{ptwise} by the one which precedes it.  Since $[L,\eta]u$ vanishes 
unless $t\le 2$ and $|x|\le 25$, the last term is clearly $O(\varepsilon)$ by \eqref{local.solution.estimate}.
The second term in \eqref{v.equation.1} is 
\begin{equation}\label{v.equation.2}
\lesssim \int_0^t \int_{|y|\ge 6(s+1)} \sum_{|\beta|+\mu\le 108} |S^\mu Z^\beta \nabla u(s,y)|
\sum_{|\beta|+\mu\le 108} |S^\mu Z^\beta \nabla^2 u(s,y)|
\:\frac{dy\:ds}{|y|}.\end{equation}

Over $|x|>6(t+1)$, we note that by (a trivial modification of) \eqref{KS.4} and \eqref{KS.6}
\begin{multline}\label{v.equation.3}
\sup_{|x|>6(t+1)} \Bigl(\la c_\kappa t-r\ra \la r\ra |P_\kappa \nabla^2 S^\mu Z^\beta u(t,x)|
+ \la r\ra \la c_\kappa t-r\ra^{1/2} |P_\kappa \nabla S^\mu Z^\beta u(t,x)|\Bigr)
\\\lesssim \sum_{|\alpha|+\nu\le 111} \|S^\nu Z^\alpha u'(t,\cd)\|_{L^2(|x|>5(t+1))}
+ \sum_{|\alpha|+\nu\le 110} \la t\ra\|S^\nu Z^\alpha Q(\nabla u,\nabla^2 u)\|_{L^2(|x|>5(t+1))}
\end{multline}
for $\kappa=1,2$ and for $|\beta|+\mu\le 108$. 
The right side of \eqref{v.equation.3} is easily seen to be $O(\varepsilon)$ from
\eqref{bdyless.solution.estimate}.

Applying this bound to the terms of \eqref{v.equation.2} and recalling that we are assuming that
$c_1=1$, it follows that \eqref{v.equation.2} is
$$\lesssim \varepsilon^2 \int_0^t \frac{1}{s^{(3/2)-}} \int_{\{y\in \ext\,:\,
|y|\ge 6(s+1)\}} \frac{1}{|y|^{3+}}\:dy\:ds \lesssim
\varepsilon^2$$
as desired.

For the second term in the left side of \eqref{v.estimate}, we use the standard energy integral
method and see that
\begin{multline*}
\partial_t \sum_{|\alpha|+\nu\le 100} \|(S^\nu Z^\alpha v)'(t,\cd)\|^2_2
+\partial_t \sum_{|\alpha|+\nu\le 100} \|\nabla\cdot(S^\nu Z^\alpha v)(t,\cd)\|_2^2
\\\lesssim \Bigl(\sum_{|\alpha|+\nu\le 100} \|S^\nu Z^\alpha v'(t,\cd)\|_2\Bigr)
\Bigl(\sum_{|\alpha|+\nu\le 100} \|S^\nu Z^\alpha Lv(t,\cd)\|_2\Bigr)
\\+\sum_{|\alpha|+\nu\le 100} \Bigl|\int_\bdy \partial_t S^\nu Z^\alpha v(t,\cd)
\nabla S^\nu Z^\alpha v(t,\cd)\cdot n\:d\sigma\Bigr|
\end{multline*}
where $n$ is the outward unit normal to $\K$ at $x\in \bdy$.  Since $\K\subset\{|x|<1\}$, we
may use \eqref{data.smallness}, \eqref{v.equation}, and a trace theorem to see that
\begin{multline}
\label{v.estimate.2.1}
\sum_{|\alpha|+\nu\le 100} \|S^\nu Z^\alpha v'(t,\cd)\|_2^2 \lesssim \varepsilon^2
\\+\Bigl(\int_0^t \sum_{|\alpha|+\nu\le 100} \|S^\nu Z^\alpha \beta\bigl(|x|/(s+1)\bigr) (1-\eta)(s,\cd)
Q(\nabla u,\nabla^2u)(s,\cd)\|_2\:ds\Bigr)^2
\\+\Bigl(\int_0^t \sum_{|\alpha|+\nu\le 100} \|S^\nu Z^\alpha [L,\eta]u(s,\cd)\|_2\:ds\Bigr)^2
\\+ \int_0^t \sum_{|\alpha|+\nu\le 101} \|S^\nu \partial^\alpha v'(s,\cd)\|^2_{L^2(|x|<1)}\:ds.
\end{multline}
Using the bound for the first term in \eqref{v.estimate}, it is easy to see that the
last term is $O(\varepsilon^2)$.  By the support properties of $\eta$ and \eqref{local.solution.estimate}, 
the preceding term satisfies the same bound.  The desired bound for the second term in the right
of \eqref{v.estimate.2.1} follows from \eqref{v.equation.3} and \eqref{bdyless.solution.estimate}
as above.

It remains only to show that the last term in the left of \eqref{v.estimate} is $O(\varepsilon)$.
Here, it suffices to show that
\begin{equation}\label{v.estimate.3.1}
\sum_{|\alpha|+\mu\le 98} \int_0^t \|\la x\ra^{-1/2} S^\mu Z^\alpha v'(s,\cd)\|^2_{L^2(|x|\le c_2s/2)}\:ds
\end{equation}
is $O(\varepsilon^2 (\log(2+t))^2)$ as the estimate follows trivially from the bound for the second term
in \eqref{v.estimate} when $|x|>c_2s/2$.

For this, we note that by \eqref{KS.5}
\begin{multline}
\label{v.estimate.3.2}
\la r\ra \la c_\kappa t-r\ra^{1/2} \sum_{|\alpha|+\mu\le 98}|P_\kappa \partial S^\mu Z^\alpha v(t,x)|\lesssim
\sum_{|\alpha|+\mu\le 100} \|S^\mu Z^\alpha v'(t,\cd)\|_2 \\+ 
\sum_{|\alpha|+\mu\le 99} (1+t)\|S^\mu Z^\alpha Lv(t,\cd)\|_2 + (1+t)\sum_{\mu\le 98} \|S^\nu v'(t,\cd)\|_{
L^\infty(|x|<2)}.
\end{multline}
The first and last terms in the right are $O(\varepsilon)$ by the bounds for the preceding terms in
\eqref{v.estimate}.  The second term in the right is
\begin{multline*}
\lesssim \sum_{|\alpha|+\mu\le 99} (1+t)\|S^\mu Z^\alpha \beta(|\cd|/(t+1)) (1-\eta)(t,\cd) Q(\nabla
u,\nabla^2 u)\|_2 \\+ \sum_{|\alpha|+\mu\le 99} (1+t)\|S^\mu Z^\alpha [L,\eta]u(t,\cd)\|_2.
\end{multline*}
Each of these terms are $O(\varepsilon)$ using \eqref{bdyless.solution.estimate} and \eqref{local.solution.estimate}
respectively.

Using the $O(\varepsilon)$ bound for the right side of \eqref{v.estimate.3.2}, it follows that
\eqref{v.estimate.3.1} is
$$\lesssim \varepsilon^2 \int_0^t \frac{1}{1+s} \|\la x\ra^{-3/2}\|^2_{L^2(|x|\le c_2s/2)}\:ds
\lesssim \varepsilon^2 (\log(2+t))^2$$
which completes the proof.\qed

The estimate \eqref{v.estimate} allows us in most instances to restrict our attention to $w-v$ which
solves
\begin{equation}
\label{wv.equation}
\begin{cases}
L(w-v)=(1-\beta)\Bigl(\frac{|x|}{t+1}\Bigr)(1-\eta)(t,x)Q(\nabla u,\nabla^2 u),\quad (t,x)\in \R_+\times\ext,\\
(w-v)|_\bdy =0,\\
(w-v)(t,x)=0,\quad t\le 0.
\end{cases}
\end{equation}
This equation meets the requirements of our estimates in the previous sections.  In particular, it has
vanishing Cauchy data and a forcing term which is supported in $|x|\le 20c_1 t$.  Depending on the estimates
being utilized, we shall use certain $L^2$ and $L^\infty$ bounds for $u$ while at other times
we shall use them for $w-v$ or $w$.  Since $u=(w-v)+v+u_0$ and since $u_0$ and $v$ satisfy 
\eqref{local.solution.estimate} and \eqref{v.estimate} respectively, it will always be the case
that a bound for $w-v$, $w$, or $u$ will imply the same bound for the others.

We now set up a continuity argument.  From this point forward, the argument resembles that of
\cite{MNS1} quite closely.
For $\varepsilon>0$ as above, we assume that we have a solution 
to \eqref{main.equation} on $0\le t\le T$ satisfying
\begin{align}
\sum_\ss{|\alpha|+\nu\le 52\\\nu\le 1} \|S^\nu Z^\alpha w'(t,\cd)\|_2 &\le A_0 \varepsilon, \label{I}\\
(1+t+r)\sum_{|\alpha|\le 40} |Z^\alpha w'(t,x)|&\le A_1 \varepsilon, \label{II}\\
(1+t+r)\sum_\ss{|\alpha|+\nu\le 55\\\nu\le 2} |S^\nu Z^\alpha (w-v)(t,x)|&\le B_1\varepsilon^2
(1+t)^{1/10} \log(2+t),\label{III}\\
\sum_{|\alpha|\le 100} \|\partial^\alpha u'(t,\cd)\|_2&\le B_2\varepsilon (1+t)^{1/40},\label{IV}\\
\sum_\ss{|\alpha|+\nu\le 72\\\nu\le 3} \|S^\nu Z^\alpha u'(t,\cd)\|_2&\le B_3 \varepsilon (1+t)^{1/20},\label{V}\\
\sum_\ss{|\alpha|+\nu\le 70\\\nu\le 3} \|\la x\ra^{-1/2} S^\nu Z^\alpha u'\|_{L^2_tL^2_x(S_t)} &\le B_4
\varepsilon (1+t)^{1/20} (\log (2+t))^{1/2}.\label{VI}
\end{align}
We may take $A_0=A_1=4C_2$ in the first two estimates, where $C_2$ is as in \eqref{v.estimate}.  For
sufficiently small $\varepsilon$, all of the above estimates hold for, say, $T=2$ by the local existence theory.
Thus, for $\varepsilon$ sufficiently small, we prove
\begin{enumerate}
\item[{\em (i.)}] \eqref{I} is valid with $A_0$ replaced by $A_0/2$,
\item[{\em (ii.)}] Assuming {\em (i.)}, \eqref{II} is valid with $A_1$ replaced by $A_1/2$,
\item[{\em (iii.)}] \eqref{III}-\eqref{VI} follow from \eqref{I} and \eqref{II} for fixed constants $B_i$.
\end{enumerate}
For such an $\varepsilon$, this then proves that a solution exists for all $t>0$.

\subsection{Preliminary results}
We begin with the following preliminary results assuming \eqref{I}-\eqref{VI}:
\begin{equation}
\label{prelim.1}
r \la c_\kappa t-r\ra |P_\kappa \partial \nabla S^\nu Z^\alpha u(t,x)|\lesssim \varepsilon (1+t)^{3/20}
\log(2+t)
\end{equation}
\begin{equation}
\label{prelim.2}
\la r\ra \la c_\kappa t-r\ra^{1/2} |P_\kappa \partial S^\nu Z^\alpha u(t,x)|\lesssim \varepsilon (1+t)^{3/20}
\log(2+t),
\end{equation}
and
\begin{equation}
\label{prelim.3}
(1+t)\|S^\nu Z^\alpha Lu(t,\cd)\|_2\lesssim \varepsilon (1+t)^{3/20} \log(2+t)
\end{equation}
for $\kappa=1,2$, $\nu\le 1$ and $|\alpha|+\nu\le 63$.

Using \eqref{KS.2} and \eqref{KS.5}, the first two estimates above follow from \eqref{prelim.3} via 
\eqref{local.solution.estimate}, \eqref{v.estimate}, \eqref{III}, and \eqref{V}.  To show the latter estimate,
we expand the left side to see that it is
$$\lesssim (1+t) \Bigl\|\sum_\ss{|\alpha|+\mu\le 32\\\mu\le 1} |S^\mu Z^\alpha u'(t,\cd)|
\sum_\ss{|\alpha|+\mu\le 64\\\mu\le 1} |S^\mu Z^\alpha u'(t,\cd)|\Bigr\|_2.$$  
By \eqref{III} and \eqref{V}, as well as \eqref{local.solution.estimate} and \eqref{v.estimate}, 
the preceding equation is clearly $O(\varepsilon^2)$.  
Notice that the same proof shows that \eqref{prelim.1} and \eqref{prelim.2}
hold with $u$ replaced by $w-v$.  

\subsection{Proof of {\em (i.)}}
As the better estimate \eqref{v.estimate} holds for $v$, it suffices to show that
\begin{equation}
\label{i.goal}
\sum_\ss{|\alpha|+\nu\le 52\\\nu\le 1} \|S^\nu Z^\alpha (w-v)'(t,\cd)\|_2^2 \lesssim \varepsilon^3.
\end{equation}

Using \eqref{e.div} with $\gamma^{IJ,jk}\equiv 0$ and with
$e_\mu[u]$ replaced by $e_\mu[S^\nu Z^\alpha (w-v)]$, we see that
the left side of \eqref{i.goal} is
\begin{multline}
\lesssim \sum_\ss{|\alpha|+\nu\le 52\\\nu\le 1} \int_0^t \int_\ext |\la \partial_t S^\nu Z^\alpha
(w-v), L S^\nu Z^\alpha (w-v)\ra| \:dy\:ds
\\+\sum_\ss{|\alpha|+\nu\le 52\\\nu\le 1} \Bigl|\int_0^t \int_\bdy
\partial_t S^\nu Z^\alpha (w-v) \partial_a S^\nu Z^\alpha (w-v) n_a \:d\sigma\:ds\Bigr|.
\end{multline}
Recalling that $\K\subset \{|x|<1\}$, and thus that the coefficients of $Z$ are $O(1)$ on $\bdy$,
it follows that the last term is bounded by
$$\sum_\ss{|\alpha|+\nu\le 53\\\nu\le 1} \int_0^t \int_{\{|x|<1\}} |S^\nu \partial^\alpha (w-v)'(s,y)|^2\:dy\:ds.$$
This boundary term is seen to be $O(\varepsilon^4)$ using \eqref{III}.

By \eqref{commutators}, \eqref{wv.equation}, and 
the fact that the vector fields preserve the null structure\footnote{See \cite{Si}, Proposition 3.1},
we may apply \eqref{null.condition.decay}, and thus it suffices to dominate
\begin{multline}
\label{energy.terms.remaining}
\int_0^t \int_\ext \frac{1}{|y|} \sum_\ss{|\alpha|+\nu\le 53\\\nu\le 1} |S^\nu Z^\alpha u|
\sum_\ss{|\alpha|+\nu\le 53\\\nu\le 1} |S^\nu Z^\alpha \nabla u|
\sum_\ss{|\alpha|+\nu\le 52\\\nu\le 1} |S^\nu Z^\alpha \partial_t (w-v)|\:dy\:ds
\\
+\int_0^t \int_\ext \sum_\ss{1\le I,J,K\le 2\\(I,J)\neq (J,K)}
\sum_\ss{|\alpha|+\nu\le 52\\\nu\le 1} |P_K \partial S^\nu Z^\alpha (w-v)|
\sum_\ss{|\alpha|+\nu\le 52\\\nu\le 1} |P_I \nabla S^\nu Z^\alpha u|\\\times
\sum_\ss{|\alpha|+\nu\le 53\\\nu\le 1} |P_J \nabla S^\nu Z^\alpha u|\:dy\:ds.
\end{multline}

For the first term in \eqref{energy.terms.remaining}, we apply \eqref{local.solution.estimate},
\eqref{v.estimate}, and \eqref{III} to see that it is
\begin{multline*}
\lesssim \varepsilon \int_0^t \frac{\log(2+s)}{(1+s)^{9/10}} \sum_\ss{|\alpha|+\nu\le 53\\\nu\le 1}
\|\la y\ra^{-1/2} S^\nu Z^\alpha u'(s,\cd)\|_{L^2_y} \\\times
\sum_\ss{|\alpha|+\nu\le 52\\\nu\le 1} \|\la y\ra^{-1/2} S^\nu Z^\alpha (w-v)'(s,\cd)\|_{L^2_y}\:ds.\end{multline*}
By the Schwarz inequality, \eqref{local.solution.estimate}, \eqref{v.estimate}, and \eqref{VI},
this term's contribution is $O(\varepsilon^3)$.

To be concrete, we will focus on the case that $I\neq K$, $I=J$ for the second term in
\eqref{energy.terms.remaining}.  A symmetric argument will then complete the proof.  For
$\delta < (c_1-c_2)/2$, 
$$\{y\in\ext\,:\, |y|\in [(1-\delta)c_Is, (1+\delta)c_Is]\}\cap 
\{y\in\ext\,:\,|y|\in [(1-\delta)c_Ks, (1+\delta)c_Ks]\}=\emptyset.$$
Thus, it suffices to show the estimate over the complements of each of these sets separately.

By \eqref{prelim.2}, it follows that over $\{|y|\not\in [(1-\delta)c_Ks, (1+\delta)c_Ks]\}$
the second term in \eqref{energy.terms.remaining} is dominated by
$$\varepsilon \int_0^t \frac{\log(2+s)}{(1+s)^{7/20}} \int_{\{|y|\not\in
[(1-\delta)c_Ks,(1+\delta)c_Ks]\}} \frac{1}{\la y\ra} \sum_\ss{|\alpha|+\mu\le 53\\\mu\le 2}
|\partial S^\mu Z^\alpha u|^2\:dy\:ds.$$
The required $O(\varepsilon^3)$ bound follows from \eqref{VI}.  A symmetric argument
then completes the proof of {\em (i.)}.

\subsection{Proof of {\em (ii.)}}
By \eqref{v.estimate}, it again suffices to show \eqref{II} when $w$ is replaced by $w-v$.  With
this substitution, it follows from \eqref{ptwise.2} that the left side of \eqref{II} is dominated
by
\begin{multline}
\label{ptwise.terms}
\sum_\ss{|\alpha|+\nu\le 46\\\nu\le 1} \int_{\theta t}^t \int_\ext |S^\nu Z^\alpha
Q(\nabla u,\nabla^2 u)|\:\frac{dy\:ds}{|y|}
\\+\sup_{0\le s\le t} (1+s)\sum_{|\alpha|\le 42} \|\partial^\alpha Q(\nabla u,\nabla^2 u)(s,\cd)\|_2.
\end{multline}

For $|y|\ge c_2s/2$, the first term is trivially $O(\varepsilon^2)$ by \eqref{I}.  Thus, it remains
to control
\begin{multline*}
\int_{\theta t}^t \int_{\{|y|\le c_2 s/2\}} \sum_\ss{|\alpha|+\nu\le 46\\\nu\le 1} |S^\nu Z^\alpha
Q(\nabla u,\nabla^2 u)|\:\frac{dy\:ds}{|y|}
\\\lesssim \sum_{\kappa_1,\kappa_2=1,2}\int_{\theta t}^t \int_{\{|y|\le c_2s/2\}}
\sum_\ss{|\beta|+\mu\le 46\\\mu\le 1} |r \la c_{\kappa_1} s-r\ra^{1/2}
P_{\kappa_1} \nabla S^\mu Z^\beta u|\\\times
\sum_\ss{|\alpha|+\nu\le 46\\\nu\le 1} |r\la c_{\kappa_2} s-r\ra P_{\kappa_2}\nabla^2 S^\nu Z^\alpha u|
\:\frac{dy\:ds}{|y|^{3} \la s\ra^{3/2}},
\end{multline*}
which is $O(\varepsilon^2)$ by \eqref{prelim.1} and \eqref{prelim.2}.

Since the second term in \eqref{ptwise.terms} is easily seen to be $O(\varepsilon^2)$ using
\eqref{I} and \eqref{II}, as well as \eqref{local.solution.estimate}, the proof of {\em (ii.)} is complete.

\subsection{Proof of {\em (iii.)}} The remainder of the proof of Theorem \ref{main.theorem} follows
from the proof given in \cite{MS1}, which we sketch for completeness.  In this section, we will 
apply the results of Section \ref{energy.estimates.section} with
$$\gamma^{IJ,jk}=-2B^{IJK}_{lmn}\partial_n u^K.$$
The necessary hypotheses \eqref{gamma.symmetry} and \eqref{gamma.smallness} follow from \eqref{nonlinear.symmetry}
and \eqref{II} respectively.  The additional hypothesis \eqref{gamma.smallness.2} also follows
from \eqref{II}.

We begin by proving \eqref{IV}.  We first examine the case $\partial^\alpha u' = \partial_t^{|\alpha|} u'$.
In the sequel, we shall use elliptic regularity to prove the general result. 

For $M\le 100$, we apply \eqref{energy.dt} and \eqref{II} to see that
\begin{equation}\label{dt.1}
\partial_t E^{1/2}_M(u)(t)\lesssim \sum_{j=0}^M \|L_\gamma \partial^j_t u(t,\cd)\|_2 + \frac{\varepsilon}{1+t}
E^{1/2}_M(u)(t).\end{equation}
For fixed $M\le 100$, we have
\begin{align*}
\sum_{j\le M} |L_\gamma \partial^j_t u|&\lesssim \sum_{|\alpha|\le 40} |\partial^\alpha u'|
\sum_{j\le M-1} |\partial_t^j \partial^2 u|
+ \sum_{|\alpha|\le M-41} |\partial^\alpha u'| \sum_{40\le |\alpha|\le M/2} |\partial^\alpha u'|\\
&\lesssim \frac{\varepsilon}{1+t} \sum_{j\le M-1} |\partial^j_t \partial^2 u|
+ \sum_{|\alpha|\le M-41} |\partial^\alpha u'| \sum_{|\alpha|\le M/2} |\partial^\alpha u'|,
\end{align*}
where, in the last line, we have applied \eqref{local.solution.estimate} and \eqref{II}.  By \eqref{ell.reg}
and a similar argument, we have
\begin{multline*}
\sum_{j\le M-1} \|\partial_t^j \partial^2 u(t,\cd)\|_2 \lesssim \sum_{j\le M} \|\partial_t^j u'(t,\cd)\|_2
+ \frac{\varepsilon}{1+t} \sum_{j\le M-1} \|\partial_t^j \partial^2 u(t,\cd)\|_2
\\+ \sum_\ss{|\alpha|\le M-41\\|\beta|\le M/2} \|\partial^\alpha u'(t,\cd) \partial^\beta u'(t,\cd)\|_2.
\end{multline*}
For $\varepsilon$ sufficiently small, the second term can be bootstrapped.  Thus, by combining the 
previous two estimates with \eqref{dt.1} and using that $\sum_{j\le M} \|\partial_t^j u'(t,\cd)\|_2\approx
E_M^{1/2}(u)(t)$ for $\varepsilon$ small, we have
\begin{equation}
\label{dt.2}
\partial_t E^{1/2}_M(u)(t)\lesssim \frac{\varepsilon}{1+t} E^{1/2}_M(u)(t)
+ \sum_\ss{|\alpha|\le M-41\\|\beta|\le M/2} \|\partial^\alpha u'(t,\cd)\partial^\beta u'(t,\cd)\|_2.
\end{equation}

For $M\le 40$, the last term drops out.  Thus, by \eqref{data.smallness} and Gronwall's inequality, we have
$$\sum_{j\le 40} \|\partial_t^j u'(t,\cd)\|_2 \lesssim \varepsilon (1+t)^{C\varepsilon},$$
which provides the same bound for
$$\sum_{|\alpha|\le 40} \|\partial^\alpha u'(t,\cd)\|_2$$
by elliptic regularity and \eqref{II}.

In order to study the case $M>40$, we shall use the following lemma.
\begin{lemma}
\label{MS.lemma}
Suppose that \eqref{local.solution.estimate}, \eqref{v.estimate}, \eqref{I}, and \eqref{II} hold.
Moreover, for $\mu\le 3$ fixed, assume that
\begin{multline}
\label{MS.hyp.1}
\sum_\ss{|\alpha|+\nu\le 100-8(\mu-1)\\\nu\le\mu-1} \|S^\nu \partial^\alpha u'(t,\cd)\|_2
+ \sum_\ss{|\alpha|+\nu\le 97-8(\mu-1)\\\nu\le\mu-1} \|\la x\ra^{-1/2} S^\nu \partial^\alpha u'\|_{L^2_tL^2_x(S_t)}
\\+\sum_\ss{|\alpha|+\nu\le 96-8(\mu-1)\\\nu\le\mu-1} \|S^\nu Z^\alpha u'(t,\cd)\|_2
+\sum_\ss{|\alpha|+\nu\le 94-8(\mu-1)\\\nu\le\mu-1} \|\la x\ra^{-1/2} S^\nu Z^\alpha u'\|_{L^2_tL^2_x(S_t)}
\\\lesssim \varepsilon (1+t)^{C\varepsilon+\sigma}
\end{multline}
and, for $M\le 100-8\mu$,
\begin{multline}
\label{MS.hyp.2}
\sum_\ss{|\alpha|+\nu\le M\\\nu\le \mu} \|S^\nu \partial^\alpha u'(t,\cd)\|_2
+ \sum_\ss{|\alpha|+\nu\le M-3\\\nu\le \mu} \|\la x\ra^{-1/2} S^\nu \partial^\alpha u'\|_{L^2_tL^2_x(S_t)}
\\+\sum_\ss{|\alpha|+\nu\le M-4\\\nu\le\mu} \|S^\nu Z^\alpha u'(t,\cd)\|_2 
+\sum_\ss{|\alpha|+\nu\le M-6\\\nu\le \mu} \|\la x\ra^{-1/2} S^\nu Z^\alpha u'\|_{L^2_tL^2_x(S_t)}
\lesssim \varepsilon (1+t)^{C\varepsilon+\sigma}
\end{multline}
for $C,\sigma>0$. Then there is a constant $C'>0$ so that
\begin{multline}
\label{MS.conclusion}
\sum_\ss{|\alpha|+\nu\le M-2\\\nu\le\mu} \|\la x\ra^{-1/2} S^\nu \partial^\alpha u'\|_{L^2_tL^2_x(S_t)}
+ \sum_\ss{|\alpha|+\nu\le M-3\\\nu\le\mu} \|S^\nu Z^\alpha u'(t,\cd)\|_2
\\+ \sum_\ss{|\alpha|+\nu\le M-5\\\nu\le\mu} \|\la x\ra^{-1/2} S^\nu Z^\alpha u'\|_{L^2_tL^2_x(S_t)}
\lesssim \varepsilon (1+t)^{C'\varepsilon+C'\sigma}.
\end{multline}
\end{lemma}

\noindent{\em Proof of Lemma \ref{MS.lemma}:}  We begin by bounding the first term 
in the left of \eqref{MS.conclusion}.  By \eqref{local.solution.estimate} and \eqref{v.estimate}, it
suffices to bound this term when $u$ is replaced by $w-v$.  By \eqref{KSS.S.d}, we thus need to
establish control of
$$\sum_\ss{|\alpha|+\nu\le M-2\\\nu\le \mu} \int_0^t \|S^\nu \partial^\alpha Lu(s,\cd)\|_2\:ds
+\sum_\ss{|\alpha|+\nu\le M-3\\\nu\le\mu} \|S^\mu \partial^\alpha Lu\|_{L^2_tL^2_x(S_t)}.$$
We will show the bound for the first term.  The estimate for the second term follows from straightforward
modifications of the argument.

For $\mu=0$, we have, as above, that
$$
\sum_{|\alpha|\le M-2} |\partial^\alpha Lu|
\lesssim \sum_{|\alpha|\le 40} |\partial^\alpha u'|\sum_{|\alpha|\le M-1} |\partial^\alpha u'|
+ \sum_{|\alpha|\le M-41} |\partial^\alpha u'| \sum_{|\alpha|\le M/2} |\partial^\alpha u'|.
$$
When $M<41$, the last term vanishes and, in this case, we have that
$$\sum_{|\alpha|\le M-2}\int_0^t \|\partial^\alpha Lu(s,\cd)\|_2\:ds
\lesssim \int_0^t \frac{\varepsilon}{1+s} \sum_{|\alpha|\le M-1} \|\partial^\alpha u'(s,\cd)\|_2\:ds.$$
The desired bound then follows from \eqref{MS.hyp.2}.  When $M\ge 41$, we apply \eqref{weighted.Sobolev}
to see that
\begin{multline*}
\sum_{|\alpha|\le M-2} \int_0^t \|\partial^\alpha Lu(s,\cd)\|_2\:ds
\lesssim \int_0^t \frac{\varepsilon}{1+s}\sum_{|\alpha|\le M-1} \|\partial^\alpha u'(s,\cd)\|_2\:ds
\\+ \sum_{|\alpha|\le \max(M-38,M/2)} \|\la x\ra^{-1/2} Z^\alpha u'\|_{L^2_tL^2_x(S_t)},
\end{multline*}
and \eqref{MS.conclusion} again follows from \eqref{MS.hyp.2}.

For $\mu>0$, we have the similar bound
\begin{multline}\label{withS}
\sum_\ss{|\alpha|+\nu\le M-2\\\nu\le \mu} |S^\nu \partial^\alpha Lu|
\lesssim \sum_{|\alpha|\le 40} |\partial^\alpha u'| \sum_\ss{|\alpha|+\nu\le M-1\\\nu\le\mu} |S^\nu
\partial^\alpha u'|
\\+ \sum_\ss{|\alpha|+\nu\le M-41\\\nu\le \mu} |S^\nu \partial^\alpha u'| \sum_{|\alpha|\le M-1} |\partial^\alpha u'|
+ \sum_\ss{|\alpha|+\nu\le M/2\\\nu\le\mu-1} |S^\nu \partial^\alpha u'|
\sum_\ss{|\alpha|+\nu\le M-1\\\nu\le\mu-1} |S^\nu \partial^\alpha u'|.
\end{multline}
We may apply \eqref{II} and
\eqref{weighted.Sobolev} to see that
\begin{multline*}
\sum_\ss{|\alpha|+\nu\le M-2\\\nu\le\mu} \int_0^t \|S^\nu \partial^\alpha Lu(s,\cd)\|_2\:ds
\lesssim \int_0^t \frac{\varepsilon}{1+s}  \sum_\ss{|\alpha|+\nu\le M-1\\\nu\le\mu}
\|S^\nu \partial^\alpha u'(s,\cd)\|_2\:ds
\\+\sum_\ss{|\alpha|+\nu\le M-39\\\nu\le\mu} \|\la x\ra^{-1/2} S^\nu Z^\alpha u'\|_{L^2_tL^2_x(S_t)}
\sum_{|\alpha|\le M-1} \|\la x\ra^{-1/2} \partial^\alpha u'\|_{L^2_tL^2_x(S_t)}
\\+ \sum_\ss{|\alpha|+\nu\le \max(M-1, M/2+2)\\\nu\le \mu-1} \|\la x\ra^{-1/2} S^\nu Z^\alpha u'\|^2_{L^2_tL^2_x(S_t)},
\end{multline*}
and the desired bounds follow from \eqref{MS.hyp.1} and \eqref{MS.hyp.2}.

Once the estimate \eqref{MS.conclusion} is established for the second term in the left, the bound
for the third term in the left follows from the arguments above using \eqref{KSS.S.Z}.

It remains to bound the second term in the left of \eqref{MS.conclusion}.  To this end, we shall apply 
\eqref{energy.S.Z}.  We first examine the $\mu=0$ case.  Notice that
\begin{align*}
\sum_{|\alpha|\le M-3} \|L_\gamma Z^\alpha u(t,\cd)\|_2
&\lesssim \sum_{|\alpha|\le 40} \|Z^\alpha u'(t,\cd)\|_\infty \sum_{|\alpha|\le M-3} \|Z^\alpha u'(t,\cd)\|_2
\\
&\qquad\qquad\qquad\qquad+ \sum_{|\alpha|\le M-43,|\beta|\le M/2} \|Z^\alpha u'(t,\cd) Z^\beta u'(t,\cd)\|_2\\
&\lesssim \frac{\varepsilon}{1+t} Y^{1/2}_{M-3,0}(t) + \sum_{|\alpha|\le \max(M-41, M/2+2)}
\|\la x\ra^{-1/2} Z^\alpha u'(t,\cd)\|^2_2,
\end{align*}
and moreover notice that the last term is unnecessary when $M<43$.  By \eqref{II}, when $M<43$, it follows
from \eqref{energy.S.Z} and Gronwall's inequality that
$$\sum_{|\alpha|\le M-3} \|Z^\alpha u'(t,\cd)\|^2_2 \lesssim (1+t)^{C\varepsilon}
\Bigl(\varepsilon^2 + \sum_{|\alpha|\le M-2} \|\partial^\alpha u'\|^2_{L^2_tL^2_x([0,t]\times\{|x|<1\})}\Bigr).$$
The desired bound then follows from that for the first term in the left of \eqref{MS.conclusion}.  
When $M\ge 43$, we similarly apply \eqref{energy.S.Z}, \eqref{II}, and Gronwall's inequality to see that
\begin{multline*}
\sum_{|\alpha|\le M-3} \|Z^\alpha u'(t,\cd)\|^2_2 \\\lesssim (1+t)^{C\varepsilon}
\Bigl(\varepsilon^2 + \sum_{|\alpha|\le \max(M-41,M/2+2)} \|\la x\ra^{-1/2} Z^\alpha u'\|^2_{L^2_tL^2_x(S_t)}
\sup_{0<s<t} Y_{M-3,0}^{1/2}(s)
\\+ \sum_{|\alpha|\le M-2} \|\partial^\alpha u'\|^2_{L^2_tL^2_x([0,t]\times\{|x|<1\})}\Bigr).
\end{multline*}
Thus, by additionally applying \eqref{MS.hyp.2} and bootstrapping the $Y^{1/2}_{M-3,0}$ term, the estimate follows.

For $\mu>0$, we argue in a fashion similar to \eqref{withS}.  Indeed, we have
\begin{multline*}
\sum_\ss{|\alpha|+\nu\le M-3\\\nu\le\mu} \|L_\gamma S^\nu Z^\alpha u(t,\cd)\|_2
\lesssim \frac{\varepsilon}{1+t} Y_{M-3-\mu,\mu}^{1/2}(t) 
\\+\sum_\ss{|\alpha|+\nu\le M-41\\\nu\le\mu} \|\la x\ra^{-1/2} S^\nu Z^\alpha u'(t,\cd)\|_2
\sum_{|\alpha|\le M-3} \|\la x\ra^{-1/2} Z^\alpha u'(t,\cd)\|_2
\\+\sum_\ss{|\alpha|+\nu\le \max(M-3,M/2+1)\\\nu\le\mu-1} \|\la x\ra^{-1/2} S^\nu Z^\alpha u'(t,\cd)\|^2_2.
\end{multline*}
Thus, by \eqref{data.smallness}, \eqref{energy.S.Z}, \eqref{II}, and Gronwall's inequality, it follows that
\begin{multline*}
\sum_\ss{|\alpha|+\nu\le M-3\\\nu\le\mu} \|S^\nu Z^\alpha u'(t,\cd)\|_2^2
\\\lesssim (1+t)^{C\varepsilon}\Bigl[\varepsilon^2 + 
\sup_{0<s<t} Y^{1/2}_{M-3-\mu,\mu}(s) \Bigl(\sum_\ss{|\alpha|+\nu\le \max(M-3,M/2+1)\\\nu\le\mu-1}
\|\la x\ra^{-1/2} S^\nu Z^\alpha u'\|^2_{L^2_tL^2_x(S_t)}
\\+ 
\sum_\ss{|\alpha|+\nu\le M-41\\\nu\le\mu} \|\la x\ra^{-1/2} S^\nu Z^\alpha u'\|_{L^2_tL^2_x(S_t)}
\sum_{|\alpha|\le M-3} \|\la x\ra^{-1/2} Z^\alpha u'\|_{L^2_tL^2_x(S_t)}\Bigr)
\\+\sum_\ss{|\alpha|+\nu\le M-2\\\nu\le\mu} \|S^\nu \partial^\alpha u'\|^2_{L^2_tL^2_x([0,t]\times\{|x|<1\})}\Bigr].
\end{multline*}
Applying \eqref{MS.hyp.1} and \eqref{MS.hyp.2} completes the proof.\qed

Reverting back to \eqref{dt.2} and applying \eqref{weighted.Sobolev} and Gronwall's inequality, it follows that
$$E_M^{1/2}(u)(t)\lesssim (1+t)^{C\varepsilon}\Bigl(\varepsilon + \sum_{|\alpha|\le \max(M-39, M/2)}
\|\la x\ra^{-1/2} Z^\alpha u'\|^2_{L^2_tL^2_x(S_t)}\Bigr)$$
for $40\le M\le 100$.  Using this in an induction argument based on the previous lemma gives 
$$\sum_{j\le 100} \|\partial_t^j u'(t,\cd)\|_2\lesssim \varepsilon (1+t)^{C\varepsilon+\sigma},$$
and coupled with elliptic regularity, this yields \eqref{IV}.  Additionally, from Lemma \ref{MS.lemma}, we
obtain 
\begin{multline}
\label{no.S}
\sum_{|\alpha|\le 98} \|\la x\ra^{-1/2} \partial^\alpha u'\|_{L^2_tL^2_x(S_t)}
+ \sum_{|\alpha|\le 97} \|Z^\alpha u'(t,\cd)\|_2
\\+\sum_{|\alpha|\le 95} \|\la x\ra^{-1/2} Z^\alpha u'\|_{L^2_tL^2_x(S_t)}\lesssim \varepsilon (1+t)^{C'\varepsilon
+C'\sigma},
\end{multline}
which gives the $\nu=0$ versions of \eqref{V} and \eqref{VI} and will provide the hypothesis
\eqref{MS.hyp.1} for further applications of the lemma corresponding to a single occurrence of $S$.

We now approach the proofs of the estimates involving the scaling vector field $S$ following a similar
strategy.  For a fixed $\mu=1,2,3$, we shall assume \eqref{MS.hyp.1}, and the main step is to show that
\begin{equation}
\label{S.energy.goal}
\sum_\ss{|\alpha|+\nu\le 100-8\mu\\\nu\le\mu} \|S^\nu \partial^\alpha u'(t,\cd)\|_2\lesssim \varepsilon
(1+t)^{\tilde{C}(\varespilon+\sigma)}
\end{equation}
for some $\tilde{C}>0$.  

We need first to establish \eqref{dt.bound}.  Noticing that, for $M\le 100-8\mu$,
\begin{multline*}
\sum_\ss{j+\nu\le M\\\nu\le\mu} \Bigl(|\tilde{S}^\nu \partial_t^j L_\gamma u| + |[\tilde{S}^\nu \partial_t^j,
\gamma^{kl}\partial_k\partial_l]u|\Bigr)
\\\lesssim \sum_{|\alpha|\le 40} |\partial^\alpha u'| 
\sum_\ss{j+\nu\le M-1\\\nu\le\mu}
|\tilde{S}^\nu\partial_t^j \partial^2 u| 
+\sum_\ss{|\alpha|+\nu\le M-40\\\nu\le\mu} |S^\nu \partial^\alpha u'| \sum_\ss{|\alpha|+\nu\le M
\\\nu\le\mu-1}
|S^\nu \partial^\alpha u'|,
\end{multline*}
it follows from elliptic regularity, \eqref{weighted.Sobolev}, and \eqref{II} that
\begin{multline*}
\sum_\ss{j+\nu\le M\\\nu\le\mu} \Bigl(\|\tilde{S}^\nu \partial_t^j L_\gamma u(t,\cd)\|_2
+ \|[\tilde{S}^\nu \partial_t^j, \gamma^{kl}\partial_k\partial_l]u(t,\cd)\|_2\Bigr)
\lesssim \frac{\varepsilon}{1+t} \sum_\ss{j+\nu\le M\\\nu\le\mu} \|\tilde{S}^\nu \partial_t^j u'(t,\cd)\|_2
\\+\sum_\ss{|\alpha|+\nu\le M-38\\\nu\le\mu} \|\la x\ra^{-1/2} S^\nu Z^\alpha u'(t,\cd)\|_2^2
+\sum_\ss{|\alpha|+\nu\le M\\\nu\le\mu-1} \|\la x\ra^{-1/2} S^\nu \partial^\alpha u'(t,\cd)\|_2^2.
\end{multline*}
Setting the sum of the last two terms above to be $H_{\nu,M-\nu}(t)$, it follows from 
\eqref{energy.L.d}, \eqref{II}, and \eqref{MS.hyp.1} that
\begin{multline}
\label{S.energy.1}
\sum_\ss{|\alpha|+\nu\le M\\\nu\le\mu} \|S^\mu \partial^\alpha u'(t,\cd)\|_2
\lesssim (1+t)^{A(\varepsilon+\sigma)}\varepsilon
\\+ (1+t)^{A\varepsilon}
\sum_\ss{|\alpha|+\nu\le M-38\\\nu\le\mu} \|\la x\ra^{-1/2} S^\nu Z^\alpha u'\|_{L^2_tL^2_x(S_t)}^2
\\+(1+t)^{A\varepsilon}
\sum_\ss{|\alpha|+\nu\le M\\\nu\le\mu-1} \int_0^t \|S^\nu \partial^\alpha u'(s,\cd)\|_{L^2(\{|x|<1\})}\:ds.
\end{multline}

To control the last term in \eqref{S.energy.1}, we apply \eqref{local.solution.estimate},
\eqref{v.estimate}, and \eqref{bdy.term} to see that it is 
$$\lesssim \varepsilon \log(2+t)(1+t)^{A\varepsilon} + (1+t)^{A\varepsilon}\sum_\ss{|\alpha|+\nu\le M+3\\\nu\le \mu-1}
\int_0^t \int |S^\nu Z^\alpha Lu(s,y)|\:\frac{dy\:ds}{|y|},$$
where we have also applied Sobolev's lemma to bound the second term in the right of \eqref{bdy.term}
by the preceding one.  Since the second term above is 
$$\le \sum_\ss{|\alpha|+\nu\le M+4\\\nu\le\mu-1} \|\la x\ra^{-1/2} S^\nu Z^\alpha u'\|^2_{L^2_tL^2_x(S_t)},$$
it follows from \eqref{MS.hyp.1} that
\begin{multline*}
\sum_\ss{|\alpha|+\nu\le M\\\nu\le\mu} \|S^\mu \partial^\alpha u'(t,\cd)\|_2
\lesssim (1+t)^{A'(\varepsilon+\sigma)}\varepsilon
\\+ (1+t)^{A\varepsilon}
\sum_\ss{|\alpha|+\nu\le M-38\\\nu\le\mu} \|\la x\ra^{-1/2} S^\nu Z^\alpha u'\|_{L^2_tL^2_x(S_t)}^2.
\end{multline*}
For $M\le 40$, this yields \eqref{S.energy.goal}.  For $M\ge 40$, we again argue inductively using
Lemma \ref{MS.lemma}.  Arguing as such proves
\begin{multline*}
\sum_\ss{|\alpha|+\nu\le 100-8\mu\\\nu\le\mu} \|S^\nu \partial^\alpha u'(t,\cd)\|_2
+ \sum_\ss{|\alpha|+\nu\le 97-8\mu\\\nu\le\mu} \|\la x\ra^{-1/2} S^\nu \partial^\alpha u'\|_{L^2_tL^2_x(S_t)}
\\+\sum_\ss{|\alpha|+\nu\le 96-8\mu\\\nu\le\mu} \|S^\nu Z^\alpha u'(t,\cd)\|_2
+\sum_\ss{|\alpha|+\nu\le 94-8\mu\\\nu\le\mu} \|\la x\ra^{-1/2} S^\nu Z^\alpha u'\|_{L^2_tL^2_x(S_t)}
\\\lesssim \varepsilon (1+t)^{C\varepsilon+\sigma}, \quad \mu=0,1,2,3,
\end{multline*}
which implies \eqref{V} and \eqref{VI}.

It remains only to prove \eqref{III}.  Using \eqref{ptwise}, it is easy to see that the left side of \eqref{III}
is dominated by the square of the left side of \eqref{VI}, from which \eqref{III} follows.  This completes
the proof of Theorem \ref{main.theorem}.

\bigskip

\end{document}